\documentclass[journal]{IEEEtran}
\usepackage[english]{babel}
\usepackage{bbm}            
\usepackage{amsmath}
\usepackage{graphicx}
\usepackage{multirow}
\usepackage{comment}
\usepackage{multicol}
\usepackage{xcolor}
\usepackage{amssymb}
\usepackage[short]{optidef}
\usepackage{hyperref}
\usepackage{url}
\usepackage{pgfplots}
\usepackage{adjustbox}

\usepackage[ruled]{algorithm2e}

\DeclareMathOperator*{\argmin}{arg\,min}
\newcommand{\vect}[1]{\boldsymbol{#1}}
\SetKwRepeat{Do}{do}{while}%

\pgfplotsset{compat=1.17}

\begin{document}

%
\title{Mixed-Integer Nonlinear Programming for State-based Non-Intrusive Load Monitoring}
%
%
%

\author{Marco Balletti, Veronica Piccialli, and~Antonio M.~Sudoso}

\IEEEoverridecommandlockouts
\IEEEpubid{\makebox[\columnwidth]{\copyright 2022 IEEE. \url{https://doi.org/10.1109/TSG.2022.3152147}}
\hspace{\columnsep}\makebox[\columnwidth]{ }}
\maketitle
\IEEEpubidadjcol
\begin{abstract}
Energy disaggregation, known in the literature as Non-Intrusive Load Monitoring (NILM), is the task of inferring the energy consumption of each appliance given the aggregate signal recorded by a single smart meter.  In this paper, we propose a novel two-stage optimization-based approach for energy disaggregation.
In the first phase, a small training set consisting of disaggregated power profiles is used to estimate the parameters and the power states by solving a mixed integer programming problem. Once the model parameters are estimated, the energy disaggregation problem is formulated as a constrained binary quadratic optimization problem. We incorporate penalty terms that exploit prior knowledge on how the disaggregated traces are generated, and appliance-specific constraints characterizing the signature of different types of appliances operating simultaneously. 
Our approach is compared with existing optimization-based algorithms both on a synthetic dataset and on three real-world datasets. The proposed formulation is computationally efficient, able to disambiguate loads with similar consumption patterns, and successfully reconstruct the signatures of known appliances despite the presence of unmetered devices, thus overcoming the main drawbacks of the optimization-based methods available in the literature.
\end{abstract}


%
\IEEEpeerreviewmaketitle

\section*{Notation}
\label{sec:notation}
Throughout the paper, we denote vectors by boldface lowercase letters and matrices by boldface uppercase letters. We also denote by $\boldsymbol{1}_n \in \mathbb{R}^{n}$ the vector of all ones of size $n$. 
Inputs:
\begin{IEEEdescription}
    \item[$N$] Number of appliances.
    \item[$S_i$] Number of power states of the appliance $i$.
    \item[$M$] Number of training time periods.
    \item[$V$] Number of validation time periods.
    \item[$T$] Number of test time periods.
    \item[$H$] Number of time intervals in $T$.
    \item[$y_i(t)$] Active power of the appliance $i$ at time $t$.
    \item[$\hat{y}_i(t)$] Estimated active power of the appliance $i$ at time $t$.
    \item[$y(t)$] Total active power of the time $t$.
    \item[$\vect{z}_i(t)$] Historical observations of the appliance $i$ up to time $t$.
\end{IEEEdescription}
Sets:
\begin{IEEEdescription}[\IEEEsetlabelwidth{$V_1,V_2,V_3,V_4,V_5$}]
    \item[$\mathcal{N} = \{1, \dots N\}$] Set of all appliances.
    \item[$\mathcal{O} \subseteq \mathcal{N}$] Set of Type-IV appliances.
    \item[$\mathcal{F} \subseteq \mathcal{N}$] Set of Type-II appliances.
    \item[$\mathcal{S}_i = \{1, \dots S_i\}$] Set of states of the appliance $i$.
    \item[$\mathcal{M} = \{1, \dots M\}$] Set of training indices.
    \item[$\mathcal{V} = \{1, \dots V\}$] Set of validation indices.
     \item[$\mathcal{T} = \{1, \dots T\}$] Set of test indices.
     \item[$\mathcal{H} = \{1, \dots ,H\}$] Set of interval indices.
     \item[$\mathcal{T}_h$] Set of indices $t \in \mathcal{T}$ related to the interval $h \in \mathcal{H}$.
\end{IEEEdescription}
Model Parameters:
\begin{IEEEdescription}[\IEEEsetlabelwidth{$V_1,V_2$}]
    \item[$\vect{p}_i$] Power levels of the appliance $i$.
    \item[$\lambda_1$] Non-negative penalty parameter for temporal smoothness.
    \item[$w_i$] Non-negative weight for the piece-wise constant consumption profile of the appliance $i$.
    \item[$\lambda_2$] Non-negative penalty parameter for device sparsity.
    \item[$l_i$] Non-negative weights inversely proportional to the ON time of the appliance $i$.
    \item[$s_i(t)$] Probability that the appliance $i$ is ON at time $t$.
    \item[\footnotesize{$a_{ij},b_{ij}$}] Minimum and maximum active time for state $j$ of the appliance $i$.
    \item [\footnotesize{$\boldsymbol{U}_i,\boldsymbol{D}_i$}] Indicator matrices defining the transition order between the states of the appliance $i$.
    \item[$o_i$] Maximum number of upward transitions for the appliance $i$.
    \item[$m_{ih}$] Maximum consumption allowed in the $h$-th time interval for the appliance $i$.
\end{IEEEdescription}
Variables:
\begin{IEEEdescription}[\IEEEsetlabelwidth{$V_1,V_2,V_3$}]
    \item [$\boldsymbol{x}_i(t)$] State indicator variable of the appliance $i$.
    \item[\footnotesize{$\boldsymbol{u}_i(t), \boldsymbol{d}_i(t)$}]Indicator variables for upward and downward transitions of the appliance $i$ at time $t$.
    \item [$\vect{\phi}_{ij}$] Autoregressive coefficients for state $j$ of the appliance $i$.
    \item [$\delta_{ij}(t)$] Assignment of $\vect{z}_i(t)$ to the $j$-th submodel of the appliance $i$.
\end{IEEEdescription}

\section{Introduction}
%
%
%
%
Non-Intrusive Load Monitoring (NILM), or energy disaggregation, is the task of estimating the energy consumption of each appliance starting from the aggregate power signal recorded by a single meter \cite{hart}. Given the aggregate power consumption and a set of target appliances, a NILM system aims to identify at each time step the active appliances and their relative contribution to the power consumption.
NILM is recognized as an important Smart Grid technology that provides energy breakdown information without the need of installing multiple monitoring devices at the appliance level.
Energy saving is arguably the most popular application of the NILM service \cite{zeifman2011nonintrusive, zoha2012non}. It can be useful for increasing energy savings for both providers and residential users.
Regarding the providers, detailed appliance usage information allow them to predict the energy demand, apply personalized management policies and service recommendations, as well as to promote future economic and environmental objectives.
Residential users, on the other hand, can obtain more awareness of how much energy being spent on their appliances and can exploit this information to take the proper actions for reducing their bills.

Since Hart’s seminal paper \cite{hart}, numerous energy disaggregation algorithms have been proposed in the literature to improve his results.
NILM approaches can be divided into event-based or state-based methods. Depending on whether the ground truth is available or not, NILM algorithms can be supervised or unsupervised \cite{he2016non}. Event-based approaches use edge detection techniques to capture statistically significant variations of the signal. Most common event-based approaches involve an unsupervised event detection of the aggregate signal and a supervised classification to assign the known appliances to the detected events. The classification tools that can be found in the literature are based on Support Vector Machines (SVM) \cite{hassan2013empirical}, Decision Trees \cite{chowdhury2019non} and a combination of various classification methods \cite{lin2014non}. In addition, clustering techniques based on Dynamic Time Warping (DTW) are used to identify windows with similar consumption patterns and to extract representative load signatures \cite{elafoudi2014power, liu2017dynamic}. 
Recent research on event-based NILM involves Graph Signal Processing (GSP) algorithms \cite{he2016non, zhao2018improving}, both supervised and unsupervised, where the spatio-temporal correlation of the data is found by embedding the signals onto a graph. Recently, multi-label classification algorithms employing time-domain and wavelet-domain features have been suggested (see \cite{tabatabaei2016toward} and reference therein).

State-based approaches, instead, consider each appliance as a finite-state machine and disaggregate the total power signal based on the appliance load distribution models. In order to explore all the possible combinations among the different appliances' state sequences, Combinatorial Optimization (CO) and Hidden Markov Models (HMM) \cite{kolter2012approximate, kong2016hierarchical, bonfigli2017non} are commonly used among the state-based approaches. However, as the number of appliances increases, the number of combinations of states increases exponentially, increasing the problem complexity. Another limitation of these approaches is that they tend to produce unsatisfactory results in presence of unknown loads.
In \cite{parson2012non} differential HMMs are used with the expectation-maximization algorithm to generate state transition models in an unsupervised manner. In \cite{makonin2015exploiting} sub-metering measurements are used to build super-state HMMs and the inference is performed  through a sparse Viterbi algorithm resulting in an efficient estimation of the energy consumption. In \cite{rahimpour2017non} a non-negative matrix factorization technique is used to decompose the aggregated signal into appliance signatures.

Among the supervised methods, in the last few years, state-of-the-art performance has been achieved by using end-to-end Deep Neural Networks (DNN) \cite{kelly2015neural, zhang2018sequence, kaselimi2020context, piccialli2021improving}. 
The main drawback of DNNs, it that they strongly rely on extended training samples of both aggregated and disaggregated data. To obtain good performance it is recommended to use measurements at the granularity of seconds which may be not available since data points at the resolution of minutes are more practical to be recorded by commercial smart meters \cite{shin2019data}. Furthermore, they focus on training one model per appliance (one-to-one mapping), resulting that a complete NILM system should integrate as many models as the number of devices the target environment contains. Thus, these types of solutions are not directly applicable in real-time situations. On the other hand, by training one model for each target appliance, DNNs are able to identify the activations and estimate the consumption of the associated appliances regardless of the number of appliances or the proportion of power consumption from non-target appliances \cite{kaselimi2020context, piccialli2021improving}.

\subsection{Related Work on Optimization}
NILM can be treated as a blind source separation problem, which tries to estimate the composition of the input from the output (one-to-many mapping). In this light, optimization-based approaches are state-based, so the studied appliances are assumed to have discrete states in their signatures \cite{mauch2015new}. 
Unlike machine learning research on NILM which mostly focuses on training one model per appliance, state-based optimization approaches can disaggregate multiple loads at the same time.
In \cite{hart}, Hart formulates the energy disaggregation problem as a combinatorial optimization (CO) problem. In his formulation, the disaggregated energy consumption is obtained by combining all the possible appliances' states so that the difference between the aggregated power and the sum of the appliance-level power is minimized. In addition to the high computational cost, the drawbacks of Hart's formulation are the possible confusion derived from loads with similar power levels, and the inability to correctly disaggregate the target appliances in presence of unknown appliances. In fact, if the measurement noise or the contribution of the unknown appliances is significant, the CO problem generates a poor solution where all the appliances are set to an active state most of the time. 

Recently, various integer programming approaches have been proposed in the literature, but very few of them have attempted to deal with the challenges of Hart's formulation. Suzuki et al. \cite{suzuki2008nonintrusive} formulated the NILM problem as an Integer Linear Programming (ILP) problem where they used current waveforms from multiple loads. Egarter and Elmenreich \cite{egarter2015load} investigated the CO approach from a theoretical point of view and discussed the equivalency with the knapsack problem. They concluded that it is hard to disaggregate loads with similar consumption patterns and proposed as future work a multi-objective optimization approach. Bhotto et al. \cite{bhotto2017} proposed several enhancements to the ILP model in \cite{suzuki2008nonintrusive}, such as always-on appliances and median filtering. Most of the enhancements in \cite{bhotto2017} were included as a pre-processing step rather than constraints. They also proposed a linear programming refinement to deal with the minimum and maximum transient spans of the power states. In \cite{piga2015sparse}, Piga et al. proposed instead a convex formulation with continuous variables and enforced sparsity by means of the Lasso penalty.
In \cite{wittmann2018}, Wittmann et al. proposed a mixed-integer programming model that exploits both active and reactive power. They added basic state machine constraints but their method is window-based and does not work on every instance of data. More recently, Zeinal-Kheiri et. al \cite{zeinal-kheiri2020} enriched the integer programming model by enforcing the power profiles of the appliances to be piece-wise constant over time. The state transitions of each appliance are modeled with non-linear constraints. Clearly, their proposed model is an integer non-linear programming problem that is difficult to solve to global optimality as the number of appliances or the scheduling horizon increase. 
In addition to deterministic optimization approaches, derivative-free search strategies such as ant colony optimization \cite{lin2013development}, differential evolution \cite{lu2019non} and genetic algorithms \cite{egarter2013evolving} have also been employed.

\subsection{Our Contribution}
According to \cite{zoha2012non}, household appliances can be categorized on the basis of their operational states: Type-I (ON/OFF), Type-II (Finite State Machines), Type-III (Continuously Variable), or Type-IV (Always-On). In this paper, we are concerned about expanding the CO problem for energy disaggregation in order to disambiguate the signatures of different types of appliances and handle the major challenges that arise due to the presence of unknown loads. Results presented in most of the optimization-based literature are obtained on an artificial signal, created by summing up the plug level power measurements of the target appliances and by subtracting the contribution of unknown appliances. This makes results not suitable for real-world applications unless nearly all the devices connected to the household main meter are known. Therefore, in this paper, special emphasis is placed on load disaggregation in an environment where the unmetered part is typically plenty and highly variable. 

We formulate the energy disaggregation problem as a least-square error minimization problem with regularization terms that promote sparsity and exploit prior knowledge on the disaggregated traces. In order to further reduce the computational complexity, we narrow down the feasible region by adding appliance-specific and signal aggregate constraints that characterize the signature of different types of appliances operating simultaneously.

The contributions of our paper are the following:
\begin{enumerate}
    \item A novel penalty-based binary quadratic programming formulation with appliance-specific constraints is
    proposed to better disambiguate the load of similar appliances in presence of unknown loads.
    \item A novel optimization-based automatic state detection algorithm is developed to estimate the power levels of the appliances and their respective transient behavior.
    \item An automatic procedure is proposed for inferring the parameters of the optimization problem by using a small training set of disaggregated data.
    \item The superiority of the approach with respect to previous works in the field of optimization-based methods for NILM is validated by extensive testing on different publicly available datasets.
    \item The deployment of the proposed algorithm in a real-world environment is discussed and evaluated experimentally.
\end{enumerate}

Our optimization-based load disaggregation algorithm consists of three steps. First, we learn the parameters of the optimization problem by using a very small training set of power consumption signatures from different devices. More in detail, the power states and their dynamics are obtained by efficiently finding a good locally optimal solution of a non-convex mixed-integer optimization problem that jointly clusters and estimates an autoregressive model for each power state. Then, taking as input the aggregated signal, we solve the binary quadratic disaggregation problem which incorporates different types of constraints and prior knowledge on the estimated consumptions. Finally, we post-process the output of the optimization problem to keep into account the transient spans of each appliance. Our method is computationally efficient for low-frequency data, i.e. 1 minute to 10 minutes granularity, which is a realistic solution for many commercially available smart meters. The remainder of this paper is organized as follows. In Section II, NILM is formulated as a binary quadratic optimization problem with linear constraints. In Section III, a mixed-integer programming formulation is described to ensure reliable automatic state detection. In Section IV, the proposed method is experimentally evaluated against the previous optimization-based methods on publicly available datasets. Finally, Section V concludes the paper.

\section{BQP Disaggregation Algorithm}
In this section, we formulate the energy disaggregation problem as an error minimization problem. Specifically, we start from Hart's formulation in \cite{hart} and we propose several enhancements aimed at improving the disaggregation accuracy. For all $t \in \mathcal{T}$, the energy disaggregation problem can be formulated as the following integer problem: 

\begin{mini}[2]
    {}{\bigg| y(t) - \sum_{i=1}^{N} \boldsymbol{p}_i^\top \boldsymbol{x}_i(t) \bigg |}
    {}{}
    \addConstraint{\boldsymbol{x}_i(t) \in \{0, 1\}^{S_i}}{\quad \forall i \in \mathcal{N}}{},
    \label{prob:original}
\end{mini}

where $\boldsymbol{p}_i \in \mathbb{R}^{S_i}$ is the vector of the non-zero Watt consumption levels of the appliance $i$ (i.e. the power states) and $\boldsymbol{x}_i(t) \in \{0, 1\}^{S_i}$ is the state variable of the appliance $i$. The $j$-th component of $\boldsymbol{x}_i(t)$ is set to 1 if the state $j$ of the appliance $i$ is active, and 0 otherwise.

The solution of the Problem (\ref{prob:original}) is generally not unique since there are many combinations of power levels resulting in the same aggregate signal. Thus, we add regularization and several constraints to better distinguish the contribution of all the appliances. 


\subsection{Objective Function}
The objective function we minimize is given by the sum of three terms. The main term is the fitting error calculated as the sum of squared differences between the aggregated power consumption and the sum of the disaggregated consumption of each appliance:

\begin{equation}
\label{eq:obj_term1}
f(\boldsymbol{x}(t)) = \sum_{t=1}^{T} \bigg( y(t) - \sum_{i=1}^{N} \boldsymbol{p}_i^\top \boldsymbol{x}_i(t) \bigg)^2.
\end{equation}

We stress that the fitting error takes into account the presence of unknown electrical loads since $y(t)$ is the aggregate power measured from the main meter. 
We add to the fitting error (\ref{eq:obj_term1}) a term to exploit the knowledge that the underlying appliance signals are piece-wise constant over time. More in detail, we enforce sparseness and temporal smoothness on the state indicator vectors $\boldsymbol{x}_i(t)$ by penalizing every change of consumption level experienced by each appliance during the optimization horizon. Penalizing the norm of the difference between two consecutive parameters $\boldsymbol{x}_i(t)$ and $\boldsymbol{x}_i(t-1)$ is commonly referred in the literature as fused Lasso \cite{tibshirani2005sparsity}. For this purpose, we use the $\ell_2$ norm to directly penalizes large changes with respect to the temporal structure:
\begin{equation}
\label{eq:obj_term2_1}
g_1(\boldsymbol{x}(t)) = \lambda_1 \sum_{t=2}^{T} \sum_{i=1}^{N} w_i \| \boldsymbol{x}_i(t) - \boldsymbol{x}_i(t-1) \|_2^2,
\end{equation}
where $\lambda_1$ is a penalty parameter. The appliance-specific weights $w_i$ are non-negative parameters used to modulate the piece-wise constant behavior  of the consumption profile. These weights are chosen depending on the appliance type and should be inversely proportional to the number of state transitions. Thus, for appliances that change state very rarely (i.e. clothes dryer and dishwasher) the weight is higher than the one of the appliances that frequently switch between states.
Finally, we add to the fitting error (\ref{eq:obj_term1}) a second penalty term that promotes robustness to noise and device sparsity by preferring configurations that use a small number of operating appliances. If the structured noise given by the contribution of the unknown appliances to the aggregated energy consumption is significant, the minimization of the quadratic term in (\ref{eq:obj_term1}) would lead appliances to be set to an active state most of the time, producing unsatisfactory results. To avoid this drawback, we impose a penalty on the number of active appliances at a given time:
\begin{equation}
\label{eq:obj_term3}
g_2(\boldsymbol{x}(t)) = \lambda_2 \sum_{t=1}^{T} \sum_{i=1}^{N} l_i \big(1 - s_i(t)\big) \|\boldsymbol{x}_i(t) \|_2^2,
\end{equation}
where $\lambda_2$ is a penalty parameter. Since $\boldsymbol{x}_i(t)$ is a binary vector with at most one component equal to 1 (see constraint (\ref{eq:constraint_at_most_one_active_state})), the sum of squared $\ell_2$ norms counts the number of devices that contribute to the aggregate power consumption. The parameter $s_i(t)$ is the probability, learned during a short intrusive period, that the appliance $i$ is active at time $t$ of the time horizon; for instance, during the night hours is unusual that sparse appliances like the dishwasher are active, so an activation during this period should be penalized much more than an activation in the usual times. On the other hand, the non-negative parameters $l_i$ are inversely proportional to the ON time of the appliance; in this way, the activation of an appliance that is often ON is penalized less than the one of an appliance that is OFF most of the time. 
The hyper-parameters $\lambda_1$ and $\lambda_2$ are tuned by the user through cross-validation for balancing the trade-off between minimizing the fitting error and maximizing the sparsity.

\subsection{Constraints}
We constrain the state variables to avoid multiple active states from the same appliance. Specifically, for each appliance we allow at most one active state at a given time:
\begin{equation}
\label{eq:constraint_at_most_one_active_state}
\boldsymbol{1}_{S_i}^\top \boldsymbol{x}_i(t) \leq 1 \quad \forall t \in \mathcal{T}, \ \forall i \in \mathcal{N} \setminus \mathcal{O}.
\end{equation}

There are some appliances that operate the whole time, i.e. are always-on. For these devices the above constraint must be true at the equality:
\begin{equation}
\label{eq:constraint_always_on}
\boldsymbol{1}_{S_i}^\top \boldsymbol{x}_i(t) = 1 \quad \forall t \in \mathcal{T}, \ \forall i \in \mathcal{O},
\end{equation}
where $\mathcal{O} \subseteq \mathcal{N}$ is the set that contains the indices of the Type-IV appliances.

Many appliances operate as finite state machines and their possible state transitions can be described by a state transition diagram. In order to explicitly model the behavior of an appliance as a finite state machine, we introduce the variables $\boldsymbol{u}_i(t) \in \{0, 1\}^{S_i}$ and $\boldsymbol{d}_i(t) \in \{0, 1\}^{S_i}$ and we add the following linear constraints:
\begin{equation}
\label{eq:constraint_link_transition}
\boldsymbol{x}_i(t) - \boldsymbol{x}_i(t-1) = \boldsymbol{u}_i(t) - \boldsymbol{d}_i(t) \quad \forall t \in \mathcal{T}, \forall i \in \mathcal{N},
\end{equation}
\begin{equation}
\label{eq:constraint_at_most_one_transition}
\boldsymbol{u}_i(t) + \boldsymbol{d}_i(t) \leq \boldsymbol{1}_{S_i} \quad \forall t \in \mathcal{T}, \forall i \in \mathcal{N},
\end{equation}
These additional variables model, respectively, an upward transition (i.e. from the OFF state to an ON state) and the downward transition (i.e. from an ON state to the OFF state). In the equation (\ref{eq:constraint_link_transition}), if the appliance $i$ changes operating mode at time $t$, the $j$-th component of the vector of the difference $\boldsymbol{x}_i(t) - \boldsymbol{x}_i(t-1)$ can be 1 or -1. In the first case, the constraint sets the $j$-th component of $\boldsymbol{u}_i(t)$ to 1 representing an upward transition, and in the latter the $j$-th component of $\boldsymbol{d}_i(t)$ to 1 representing a downward transition. Constraint (\ref{eq:constraint_at_most_one_transition}) prevents the $j$-th component of both $\boldsymbol{u}_i(t)$ and $\boldsymbol{d}_i(t)$ to be simultaneously 1.

We expect an appliance to stay in a state for at least a few instants:
\begin{equation}
\label{eq:constraint_min_in_state}
\sum_{\tau=t}^{t+a_{ij}-1}  x_{ij}(\tau) \geq a_{ij} u_{ij}(t) \quad \forall t \in \mathcal{T}, \forall i \in \mathcal{N}, \forall j \in \mathcal{S}_i,
\end{equation}
where $a_{ij} \in \mathbb{R}$ is a parameter representing the minimum active time for each state $j$ of the appliance $i$. Constraint (\ref{eq:constraint_min_in_state}) enforces the variables $x_{ij}(\tau)$ to be 1 for at least $a_{ij}$ epochs when the appliance $i$ goes into state $j$ at time $t$, i.e, $u_{ij}(t) = 1$. On the other hand many appliances do not stay in the same state for a long time so we add the following constraint to force the appliances to change state before the maximum active time for that state:
\begin{equation}
\label{eq:constraint_max_in_state}
\sum_{\tau=t}^{t+b_{ij}} x_{ij}(\tau) \leq b_{ij} \quad \forall t \in \mathcal{T}, \forall i \in \mathcal{N}, \forall j \in \mathcal{S}_i,
\end{equation}
where $b_{ij} \in \mathbb{R}$ is a parameter representing the maximum active time for state $j$ of the appliance $i$.
By constraining the appliance usage with the minimum and the maximum duration we can better disambiguate appliances with similar consumption profiles.
For Type-II appliances, we can easily model a given state being active only if another state of the same appliance has finished:

\begin{equation}
\label{eq:constraint_state_machine}
\boldsymbol{U}_i\boldsymbol{u}_i(t) = \boldsymbol{D}_i\boldsymbol{d}_i(t) \qquad \forall t \in \mathcal{T}, i \in \mathcal{F},
\end{equation}
where $\mathcal{F} \subseteq \mathcal{N}$ is the set that contains the indices of the Type-II appliances, $\boldsymbol{U}_i \in \{0, 1\}^{S_i \times S_i}$ and $\boldsymbol{D}_i \in \{0, 1\}^{S_i \times S_i}$ are indicator matrices defining the transition order between the states of the appliance $i$. If there is a transition from state $j_1$ to $j_2$, then the $j_2$-th component of $\boldsymbol{u}_i(t)$ is set to 1 only if the $j_1$-th component of $\boldsymbol{d}_i(t)$ is set to 1, meaning that the previous state of the appliance $i$ is not active anymore.

Type-III appliances are the most challenging loads to disaggregate, as they arbitrarily change their power consumption \cite{wichakool2014smart}. Since we can not fully characterize Type-III appliances through explicit constraints, we treat them as ON/OFF multi-state appliances.

In order to limit the presence of spurious activations, we impose an upper bound $o_i$ on the number of upward transitions for each appliance: 
\begin{equation}
\label{eq:constraint_upward}
\sum_{t=2}^{T} \vect{1}_{S_i}^\top \vect{u}_i(t) \leq o_i \quad \forall i \in \mathcal{N}.
\end{equation}
We notice that constraint (\ref{eq:constraint_upward}) could be enforced by either using the upward transitions or the downward transitions. In fact, from the constraint (\ref{eq:constraint_link_transition}) follows that the number of upward and downward transitions in the scheduled horizon differ at most by 1 in absolute value.

The following signal-aggregate constraint imposes that the sum of the disaggregated consumption profiles does not exceed the total power measured by the main meter:


\begin{equation}
\label{eq:constraint_total_consumption}
\sum_{i=1}^{N }\boldsymbol{p}_i^\top \boldsymbol{x}_i(t) \leq y(t) \qquad \forall t \in \mathcal{T}.
\end{equation}

To better characterize the different consumption patterns in a specific period (e.g., typically used appliances in the daytime and are unlikely to be used in the nighttime and vice-versa), we impose appliance-specific constraints on the maximum power consumption allowed in the scheduled horizon. Rather than making only one constraint for each appliance in the interval $[1, T]$, a series of tighter constraints could be enforced by partitioning the scheduling horizon into a fixed number of time intervals.
As a result, the following constraint imposes that the energy consumption in a time period is less or equal than the maximum consumption $m_{ih}$ allowed in that period:
\begin{equation}
\label{eq:constraint_max_consumption_zone}
\sum_{t \in \mathcal{T}_h} \boldsymbol{p}_i^\top \boldsymbol{x}_i(t) \leq m_{ih} \qquad \forall h \in \mathcal{H}, \small \forall i \in \mathcal{N},
\end{equation}
where $\mathcal{H} = \{1, \dots, H\}$ is the set of $H$ time intervals, and $\mathcal{T}_h$ is the set of indices $t \in \mathcal{T}$ related to the $h$-th time interval. Increasing $H$ leads to tighter constraints but also implies a large number of upper bounds to estimate. Furthermore, when using limited historical observations, large $H$ may lead to $m_{ih}$ estimates varying too much in time so that the resulting accuracy may be compromised due to overfitting on the training set. In the following, we simply consider two time slots: from 1 AM to 5 AM for the nighttime and from 6 AM to 12 PM for the daytime.

In order to reduce the computational cost of the algorithm, we directly set to zero some variables by analyzing the aggregate consumption and the power levels.  More in detail, let $\mathcal{Z} = \{t \in \mathcal{T}, i \in \mathcal{N}, j \in \mathcal{S}_i \mid y(t) \ll p_{ij} \}$ be the set of indices of the appliances that are more likely to be turned off, then we simply set $x_{ij}(t) = 0$, for all $(t, i, j) \in \mathcal{Z}$. 

The overall optimization-based algorithm can be expressed as:

\begin{equation}
\label{prob:full_problem}
  \begin{split}
    \text{min }  & \sum_{t=1}^{T} \bigg( y(t) - \sum_{i=1}^{N} \boldsymbol{p}_i^\top \boldsymbol{x}_i(t) \bigg)^2 \\
    & + \lambda_1 \sum_{t=2}^{T} \sum_{i=1}^{N} w_i \|\boldsymbol{x}_i(t) - \boldsymbol{x}_i(t-1)\|_2^2 \\
    & + \lambda_2 \sum_{t=1}^{T} \sum_{i=1}^{N} l_i\big(1 - s_i(t)\big)\|\boldsymbol{x}_i(t) \|_2^2 \\ \\
    \text{s.t. } &
    \begin{aligned}[t]
      & \textrm{Constraints  (\ref{eq:constraint_at_most_one_active_state}) - (\ref{eq:constraint_max_consumption_zone})} \\
      & \boldsymbol{x}_i(t), \ \boldsymbol{u}_i(t), \ \boldsymbol{d}_i(t) \in \{0, 1\}^{S_i} & \forall t \in \mathcal{T}, \ \forall i \in \mathcal{N}.
    \end{aligned}
  \end{split}
\end{equation}

Problem (\ref{prob:full_problem}) is an optimization problem with binary variables, quadratic objective function and linear constraints. The binary quadratic program with linear constraints (BQP) is a general class of optimization problems that are known to be very difficult to solve due to the non-convexity and the integrality of the variables. The number of variables and constraints is $O(T\sum_{i=1}^{N} S_i)$. Due to its combinatorial nature, the problem becomes expensive to solve as the scheduling horizon $T$ and the number of appliances $N$ increase. However, this formulation is computationally efficient for low-frequency data where the resolution ranges from 1 minute to 10 minutes. 
Smart meter data with this granularity are more common in practical applications due to the relatively low hardware cost. In this case, Problem (\ref{prob:full_problem}) can be globally solved in order of seconds by state-of-the-art integer programming solvers without specialized hardware. 


\section{State detection and parameter estimation}
\label{sec:parameter_estimation}
The optimization model in (\ref{prob:full_problem}) relies on several input parameters that need to be estimated. The most important ones are the power states $\vect{p}_i$, representing the steady-state ratings.
Steady-state ratings are usually extracted from data sheets or by means of a clustering algorithm on the ground truth data \cite{wittmann2018, zeinal-kheiri2020}, and used to characterize the typical consumption signatures by a constant value representing the average power. 
In practice, the power states of household appliances fluctuate within a range, so that approximating them with a straight line may affect the quality of the disaggregation. 
As observed in \cite{bhotto2017}, being able to capture appliances' dynamics in addition to the steady-state ratings could be effective to increase the disaggregation accuracy. However, incorporating these dynamics as additional states in the optimization problem would dramatically increase the computational cost.
For this reason, we develop a novel clustering-based approach for steady-state ratings and transient spans extraction using a very small training set of disaggregated data. This is done by solving a mixed-integer optimization problem that jointly classifies the power levels into $S_i$ clusters and estimates an autoregressive (AR) submodel of order $q_i$ for each cluster. The autoregressive coefficient are incorporated in the post-processing phase to refine the flat estimate of each appliance. Formally, for all $i \in \mathcal{N}$, we solve:


\begin{mini}[2]
    {}{ \sum_{j=1}^{S_i}\sum_{t=q_i+1}^{M} \delta_{ij}(t) \big(y_i(t) - \boldsymbol{\phi}_{ij}^{\top} \boldsymbol{z}_i(t) \big)^2 }
    {}{}
    \label{prob:state_detection}
    \addConstraint{\vect{1}_{S_i}^\top \vect{\delta}_{i}(t)}{= 1}{\quad \forall t \in \mathcal{M}}
    \addConstraint{\delta_{ij}(t) \in \{0, 1\}}{}{\quad \forall j \in \mathcal{S}_i, \forall t \in \mathcal{M}}
    \addConstraint{\phi_{ij}^h \in \mathbb{R}}{}{\quad \forall j \in \mathcal{S}_i, \forall h \in \{0, \dots, q_i\}},
\end{mini}

where $\boldsymbol{z}_i(t) = [1, y_i(t-1), \dots , y_i(t-q_i)]^\top$ is the input vector of the historical observations and $\boldsymbol{\phi}_{ij} = [\phi^0_{ij}, \phi^1_{ij}, \dots , \phi^{q_{i}}_{ij}]^\top$ is the vector of the autoregressive coefficients.
Each binary variable $\delta_{ij}(t)$ decides whether the data point $\vect{z}_i(t)$ is assigned to the $j$-th submodel of the appliance $i$, under the constraint that each data point must be assigned to only one submodel. 
Problem (\ref{prob:state_detection}) is a mixed-integer non-linear program that can be shown to be equivalent to its continuous relaxation, where the binary constraint on $\delta_{ij}$ is replaced by $\delta_{ij} \geq 0$. This equivalence derives from the optimal solution of the relaxation being integer. The resulting continuous problem is still non-convex. In order to efficiently produce a good local optimum, we use an iterative two-block Gauss-Seidel decomposition method \cite{grippo2000convergence}. In fact, the structure of the objective function and the constraints imply that the computation of the global minimum with respect to each block can be done in a computationally efficient way. More in detail, by fixing the variables $\delta_{ij}(t)$, Problem (\ref{prob:state_detection}) becomes the least-squares estimation of $S_i$ autoregressive models of order $q_i$:

\begin{align}
\label{prob:decomp_phi}
    \boldsymbol{\phi}_{ij}^*
    & = \argmin_{\boldsymbol{\phi}_{ij}} \bigg[ \sum_{j=1}^{S_i}\sum_{t=q_i+1}^{M} \overline{\delta}_{ij}(t) \big(y_i(t) - \boldsymbol{\phi}_{ij}^{\top} \boldsymbol{z}_i(t) \big)^2 \bigg].
\end{align}

Problem (\ref{prob:decomp_phi}) is an unconstrained convex optimization problem that can be solved in closed form by setting the gradient equal to zero. On the other hand, given $\vect{\phi}_{ij}$, the problem with respect to $\delta_{ij}(t)$ is linear programming problem that is separable into $M - (q_i + 1)$ subproblems. In particular, the $t$-th subproblem has the form:

\begin{mini}[2]
    {}{ \sum_{j=1}^{S_i} \delta_{ij}(t) \big(y_i(t) - \overline{\vect{\phi}}_{ij}^{\top} \boldsymbol{z}_i(t) \big)^2 }
    {}{}
    \addConstraint{\vect{1}_{S_i}^\top \vect{\delta}_{i}(t)}{= 1}{}
    \addConstraint{\delta_{ij}(t) \geq 0}{}{\quad \forall j \in \mathcal{S}_i}.
    \label{prob:decomp_delta}
\end{mini}

The optimal solution of (\ref{prob:decomp_delta}) can be constructed by observing that if $y_i(t)$ is close to the power level $j^*$, then ${\delta}_{ij^*}(t)$ is set to 1, whereas the variables with $j \neq j^*$ are set to 0.
The algorithm starts with an initial guess of the power levels $\vect{p}_i^0$ and minimizes over one block of variables with the other fixed, and vice versa until convergence. The convergence of this alternating minimization procedure can be proved by adapting the $k$-means convergence proof that can be found in \cite{kmeans-conv}. The main idea is that the variables $\delta_{ij}(t)$ can only assume a finite number of values (i.e., the number of possible assignments is finite). Therefore, since the algorithm finds the global minimum with respect to $\vect{\phi}_{ij}$ for a given assignment, there will exist a subsequence of iterations where power levels and assignment variables do not change, satisfying the stopping criterion $\|\overline{\boldsymbol{p}}_i^k - \overline{\boldsymbol{p}}_i^{k-1}\|_2 \leq \epsilon$.

An initial estimate of the power levels can be automatically obtained by a clustering algorithm, such as $k$-means \cite{lloyd1982least, piccialli2021sos}. Since the number of power levels is unknown, $k$-means is run with an increasing number of clusters. The number of power levels is then validated with the elbow method and the cluster centers, representing the average power consumption, are selected as the initial guess $\vect{p}_i^0$ of the operating modes. The overall decomposition method for automatic state detection is illustrated in Algorithm \ref{alg:state_detection}. 

In the initialization phase, the parameters $\boldsymbol{\phi}_{ij}$ of the autoregressive subprocess are unknown, so the binary variables $\delta_{ij}(t)$ are chosen so that the squared distance between the power demanded by the appliance $i$ at time $t$ and the power levels $\vect{p}_i^0$ is minimized. In the main loop, the AR coefficients are estimated and the variables $\delta_{ij}(t)$ are adjusted on the basis of the current parameters.
Finally, for each appliance $i$ we refine the estimate of the power level $j$ by simply taking the conditional expectation of the $j$-th autoregressive subprocess.
The algorithm terminates when the difference between the power levels of two consecutive iterations is less or equal than a small value $\varepsilon$ and returns the estimated power levels $\vect{p}_i$ and the AR coefficients $\vect{\phi}_{ij}$.

\begin{algorithm}
\footnotesize
\label{alg:state_detection}
\KwData{The training set $y_i(t)$, and the tolerance $\varepsilon = 10^{-5}$.}
\KwResult{The estimated power  $\overline{\boldsymbol{p}}_i$ and the AR coefficients $\boldsymbol{\phi}_{ij}$.}
$\overline{\vect{p}}_i^0 \leftarrow \textrm{run the k-means algorithm on $y_i(t)$}$\\
\For{$t = q_i+1 \dots M$}{
    $j^* \leftarrow \argmin_{j=1 \dots S_i} \big(y_i(t) - \overline{p}_{ij}^0\big)^2$\\
    $\delta_{i, j}^0(t) = \begin{cases}
        1   & \text{if } j = j^*\\
        0   & \text{otherwise}
    \end{cases}$\
}
$k \leftarrow 1$\\

\DontPrintSemicolon

\Do{$\|\overline{\boldsymbol{p}}_i^k - \overline{\boldsymbol{p}}_i^{k-1}\|_2 > \epsilon$}{
    $\boldsymbol{\phi}_{ij}^k \leftarrow \argmin \sum_{j=1}^{S_i}\sum_{t=q_i+1}^{M} \delta_{i, j}^{k-1}(t) \big(y_i(t) - \boldsymbol{\phi}_{ij}^{\top} \boldsymbol{z}_i(t) \big)^2$\
    
    \For{$t = q_i+1 \dots M$}{
        $j^* \leftarrow \argmin_{j=1 \dots S_i} \big(y_i(t) - \boldsymbol{\phi}_{ij}^{k \top} \boldsymbol{z}_i(t)\big)^2$\\
        $\delta_{ij}^k(t) = \begin{cases}
            1   & \text{if } j = j^*\\
            0   & \text{otherwise}
        \end{cases}$\\
    }
    
    \For{$j = 1 \dots S_i$}{
        $\overline{p}_{ij}^k \leftarrow \frac{\phi_{ij}^{k, 0}}{1 - \sum_{h=1}^{q_i} \phi_{ij}^{k, h}}$\\
    }
    
    $k \leftarrow k + 1$\\
}

\KwRet{$\overline{\boldsymbol{p}}_i \text{ and } \boldsymbol{\phi}_{ij}$}
\caption{Automatic State Detection}
\end{algorithm}

Given the optimal solution of (\ref{prob:full_problem}) $\boldsymbol{x}_i^{*}(t)$, the disaggregated power consumption of the appliance $i$ at time $t$ is given by $\hat{y}_i(t) = \boldsymbol{p}_i^\top \boldsymbol{x}_i^{*}(t)$. This estimate relies only on the static power levels. In order to refine the estimate we need to post-process the output by including the transient spans through the AR coefficients. This is done by applying a moving horizon forecasting procedure for $t=q+1, \dots, T$: 
\begin{equation}\label{eq:postp}
    \hat{y}_i(t) = \max\bigg\{0,  \sum_{j=1}^{S_i} x_{ij}^{*}(t) \big(\phi_{ij}^0 + \sum_{k=1}^{q_i} \phi_{ij}^k \hat{y}_i(t-k)\big)\bigg\}.
\end{equation}

Similarly to the power states, the remaining parameters in (\ref{prob:full_problem}) can be estimated from time series analysis and limited historical information. 
The extraction of some of them may be done by visual inspection. However, automatic extraction becomes necessary when applying the model to a real-world scenario. Furthermore, for practical NILM applications, the training phase should be as short and simple as possible \cite{kim2011unsupervised}.
To this end, for $o_i, m_{ih}, a_{ij}, b_{ij}$ we compute the frequency distribution over the training set and we obtain the estimate by selecting a suitable percentile.
In particular, the parameters $o_i$ and $m_{ih}$ are computed for each interval of length $T$ over the training set and estimated by taking the 95th percentile. Likewise, the minimum and maximum time in state $a_{ij}$ and $b_{ij}$ are estimated by looking at the left and the right tails of the histogram and taking, respectively, the 5th and 95th percentile of the time spent in state $j$ in the interval $T$. The parameters $w_i$ are inversely proportional to the frequency of state transitions that characterize the device, therefore they are calculated as follows:
\begin{equation}
\label{eq:param_w}
    w_i = \frac{M}{r_i} \qquad \forall i \in \mathcal{N},
\end{equation}
where $r_i$ is the number of state transitions in the training set for the appliance $i$.
The parameters $s_i(t)$ represent the probability that the device $i$ is active (with respect to a threshold of 10 Watt) at time $t$ of the scheduling horizon. These probabilities can be estimated as follows:
\begin{equation}
\label{eq:param_s}
    s_i(t) = \frac{v_i(t)}{M/T} \qquad \forall t \in \mathcal{T}, i \in \mathcal{N},
\end{equation}
where $v_i(t)$ is the number of times that the appliance $i$ is ON at time $t$, whereas the denominator represents the number of time horizons in the training set.
Finally, the parameters $l_i$, which are inversely proportional to how much the device is typically used, are calculated as:
\begin{equation}
\label{eq:param_l}
    l_i = \frac{M}{\sum_{t=1}^M v_i(t)} \qquad \forall i \in \mathcal{N}.
\end{equation}


The validation set is used only for selecting and tuning the regularization parameters $\lambda_1 \in \Lambda$ and $\lambda_2 \in \Lambda$ in (\ref{prob:full_problem}) by means of a grid-search procedure over the grid $\Lambda \times \Lambda$. Also in the validation phase, the sub-metered power consumptions are assumed to be available.

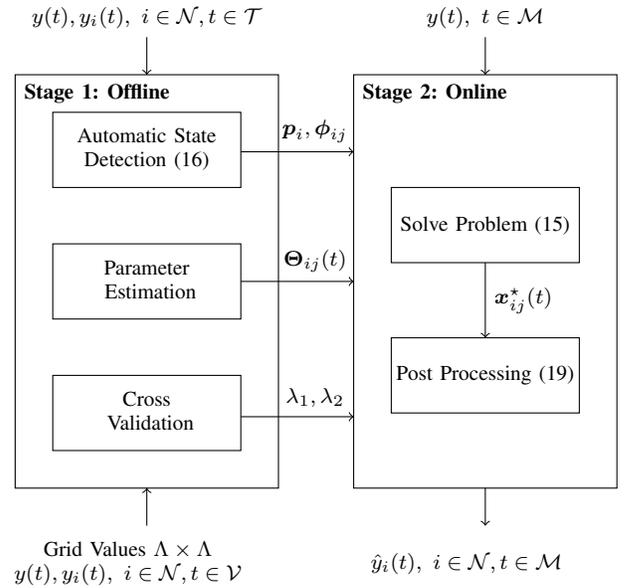
\begin{figure}
\begin{tikzpicture}[scale=1.0]
 \draw[draw=none] (1,1.5) rectangle (3.5,1)node[midway]{\footnotesize $\begin{array}{c}y(t),y_i(t),\ i \in\mathcal{N},t\in \mathcal{T}\end{array} $};
 \draw[draw=none] (5.5,1.5) rectangle (8,1)node[midway]{\footnotesize $\begin{array}{c}y(t),\ t\in {\cal M}\end{array} $};
 \draw[->](6.75,1)--(6.75,0.5);
 \draw (1,0) rectangle (3.5,-1)node[midway]{\footnotesize $\begin{array}{c}\mbox{Automatic State}\\ \mbox{Detection (\ref{prob:state_detection})} \end{array} $};
 \draw (3.5,-2.75) rectangle (1,-1.75)node[midway]{\footnotesize $\begin{array}{c}\mbox{Parameter}\\ \mbox{ Estimation} \end{array} $};
 \draw (1,-3.50) rectangle (3.5,-4.50)node[midway]{\footnotesize $\begin{array}{c}\mbox{Cross}\\ \mbox{ Validation} \end{array} $};
 \draw[draw=none] (5.5,-5.5) rectangle (7.5,-6.5)node[midway]{\footnotesize $\hat{y}_{i}(t),\ i\in \mathcal{N}, t\in \mathcal{M}$ };
  \draw[draw=none] (1,-5.5) rectangle (3,-6.5)node[midway] {\footnotesize$\begin{array}{c}\mbox{Grid Values } \Lambda\times\Lambda\\\mbox{\footnotesize $\begin{array}{c}y(t), y_i(t), \ i \in\mathcal{N},t\in \mathcal{V}\end{array} $}\end{array}$};
  \draw[->](2.25,-5.5)--(2.25,-5);
  \draw (5.5,-1) rectangle (8,-2) node[midway]{{\footnotesize Solve Problem \eqref{prob:full_problem}}};
  \draw[->] (6.75,-2)--node[right]{\footnotesize{$\boldsymbol{x}^{\star}_{ij}(t)$}}(6.75,-3);
   \draw (5.5,-3) rectangle (8,-4) node[midway]{{\footnotesize Post Processing \eqref{eq:postp}}};
\draw[->](6.75,-5)--(6.75,-5.5);
\draw[->](2.25,1)--(2.25,0.5);


\draw[->](3.5,-0.525)--node[above]{\ \quad\footnotesize $\boldsymbol{p}_i, \boldsymbol{\phi}_{ij}$}(5,-0.525);

\draw[->](3.5,-2.25)--node[above]{\ \quad \footnotesize $\boldsymbol{\Theta}_{ij}(t)$}(5,-2.25);

\draw[->](3.5,-4.05)--node[above]{\ \quad\footnotesize $\lambda_1, \lambda_2$}(5,-4.05);

\draw (0.5,0.5)node[below right] {\footnotesize{\bf Stage 1: Offline}} rectangle (4,-5);
\draw (5,0.5) node[below right] {\footnotesize{\bf Stage 2: Online}}rectangle (8.5,-5);\end{tikzpicture}
\caption{Overall disaggregation process: Stage 1 (Offline) and Stage 2 (Online). For all $i,j,t$ we denote by $\vect{\Theta}_{ij}(t) = [w_i, l_i, s_i(t), a_{ij}, b_{ij}, o_i, m_{ih}]^\top$ the vector of the parameters.}
 \label{fig:flowchart}
\end{figure}

The overall NILM algorithm can be summarized by the diagram reported in Figure \ref{fig:flowchart}.

\section{Results}
In this section we present the experiments carried out to test the performance of the overall NILM algorithm and we address several research questions:
\begin{enumerate}
    \item How does the proposed algorithm compare to state-based approaches available in the literature?
    \item Is the approach robust to the presence of unmetered appliances?
    \item How reliable is the parameter estimation procedure and how often does one need to update the parameters to reflect changes in the system?
\end{enumerate}
Overall, by answering these research questions we aim to provide insights on how the proposed algorithm can be used in practice.
We compare our NILM-BQP with NILM algorithms proposed for low sampling rates and active power measurements. Specifically, we compare with Hart's formulation (CO) \cite{hart}, Aided Linear Programming (ALIP) \cite{bhotto2017}, Sparse Optimization (SO) \cite{piga2015sparse}, State Transition Integer Programming (STIP) \cite{zeinal-kheiri2020} and Super State HMM (SHMM) \cite{makonin2015exploiting}. We first evaluate the robustness to noise on a synthetic dataset to assess how the presence of unmetered appliances affects the quality of the disaggregation, and then we evaluate our algorithm on 5 houses coming from 3 real-world datasets. Finally, we evaluate the reliability of the parameter estimation procedure and how the freshness and staleness of the parameters impact the accuracy of the disaggregation.

\subsection{Performance Metrics}
In order to evaluate the performance of our method, we use regression and classification metrics. Regression metrics measure how closely the energy consumption of an appliance matches the energy consumption predicted for that appliance, whereas the classification metrics measure how accurately NILM algorithms can predict what appliance is running in each state. Following \cite{makonin2015nonintrusive}, we use the Estimation Accuracy (EA) for the appliance $i \in \mathcal{N}$, and the Overall Estimation Accuracy (OEA) for $N$ appliances:

\begin{align}
    EA_i & = 1 - \frac{\sum_{t=1}^T |y_i(t) - \hat{y}_i(t) |}{2 \sum_{t=1}^T y_i(t)}, \\
    OEA & = 1 - \frac{\sum_{i=1}^N \sum_{t=1}^T |y_i(t) - \hat{y}_i(t) |}{2 \sum_{i=1}^N \sum_{t=1}^T y_i(t)}. 
\end{align}

In order to calculate the accuracies of non-binary classifications, we use the Finite State F-score $(FS_i)$ that is the harmonic mean of the adjusted precision $(P_i)$ and recall $(R_i)$ for the appliance $i$:

\begin{align}
    P_i = \frac{TP_i - A_i}{TP_i + FP_i}, \quad R_i = \frac{TP_i - A_i}{TP_i + FN_i}, \\
    A_i = \sum_{t=1}^{T} \frac{ | \hat{s}_i(t) - s_i(t) |}{S_i}, \quad FS_i = 2 \frac{P_i \cdot R_i}{P_i + R_i},
\end{align}

where $TP_i$, $FP_i$, and $FN_i$ stand for true positives, false positives, and false negatives for the appliance $i$, respectively. 
In the NILM context, $TP_i$ counts how many times the algorithm correctly predicts that the appliance is ON, $FP_i$ how many times the appliance is OFF but predicted ON, and $FN_i$ how many times the appliance is ON but predicted OFF.
The quantity $A_i$ is the inaccuracy portion of true positives which converts them into a discrete measure from a binary one, $\hat{s}_i(t)$ is the estimated state of the appliance $i$ at time $t$, $s_i(t)$ is the ground truth state, and $S_i$ is the number of states for the appliance $i$. The Overall FS-Fscore (OFS) is obtained by summation over all the appliances for each $TP_i$, $FP_i$, $FN_i$, and $A_i$.

To measure the contribution of unmetered loads, we report the percentage of noise in each test. The percent-noisy measure \cite{makonin2015nonintrusive} $(\%-NM)$ can be calculated on the ground truth data as follows:
\begin{equation}\label{eq:noise}
    \%-NM = \frac{\sum_{t=1}^{T} | y(t) - \sum_{i=1}^{N} y_i(t) |}{\sum_{t=1}^{T} y(t)}.
\end{equation}

\subsection{Experimental Setup}
The proposed algorithm assumes low-sampling rate measurements. We use observations at 60 seconds resolution, and we downsample the data when needed to satisfy this assumption. We use the same autoregressive order $q_i = 3$ for each appliance, even though the estimation of the transient spans could be further improved by tuning this parameter. We adopt the same experimental design for all the methods: we use only two weeks for the estimation of the parameters, a week to cross-validate the choice of the hyperparameters $\lambda_1$ and $\lambda_2$, and a week as test set. The scheduling horizon is set to 1 day ($T=1400$) which is suitable for decision making using low-frequency data \cite{shin2019data}. Furthermore, in order to minimize any bias effect, we run all the NILM algorithms one day at a time to obtain the disaggregated traces for the test week (i.e., we solve 7 optimization problems and we average the results). 
The values of $\lambda_1$ and $\lambda_2$ are determined offline through a grid-search procedure by maximizing the sum of OEA and OFS on the validation set. The grid values used in the experiments are in $\Lambda = \{200, 300, \dots, 2200\}$ and their optimal values for each dataset are reported in Table \ref{tab:lambda}.
BQP is implemented in AMPL \cite{fourer1990modeling} and solved to global optimality with Gurobi optimizer \cite{gurobi}. All the experiments have been performed on a laptop with Intel Core i7-8565U CPU and 8 GB of RAM. The source code and the parameters of all the considered case studies are available at \url{https://github.com/antoniosudoso/nilm-bqp}.

\begin{table}
\centering
\caption{Optimal hyper-parameters $\lambda_1$ and $\lambda_2$}
\scalebox{0.99}{
\begin{tabular}{|c|c|c|}
        \hline
        Dataset & $\lambda_1$ & $\lambda_2$ \\\hline\hline
        AMPDS & 1000 & 2000\\\hline
        UKDALE 1 & 500 & 1300\\\hline
        UKDALE 2 & 800 & 1500\\\hline
        REFIT 3 & 800 & 1600\\\hline
        REFIT 9 & 300 & 900\\\hline
    \end{tabular}}
\label{tab:lambda}
\end{table}

\subsection{Robustness to Noise}
It is well known that the presence of unknown appliances negatively affects the disaggregation accuracy of known appliances. Therefore a critical aspect for the deployment of the NILM service in a real-world scenario is the robustness to noise. In order to test the effectiveness of our formulation when the noise increases, we use the artificial dataset SYND \cite{klemenjak2020synthetic}. It contains 180 days of data at the granularity of 200 milliseconds of a single household with 21 appliances. This dataset has been generated using appliance-specific signatures of real-world NILM datasets. We generate artificial aggregates with an increasing number of appliances. Specifically, we disaggregate the 4 top-consuming appliances that are dishwasher, electric stove, washing machine, and iron.  
Starting from the ideal case where no noise is present, i.e., the aggregate signal is given by the sum of the 4 top-consuming appliances, we add the traces of the remaining devices one at a time in decreasing order of energy consumption, and we treat them as unknown appliances. In this way, we generate 18 different datasets, corresponding to artificial aggregates of increasing noise scenarios, where the percentage of noise ranges from $0\%$ to around $30\%$. The noisiest dataset contains all the remaining 17 appliances treated as noise, resulting in the original artificial aggregate signal. 

\begin{figure}
\centering
\begin{tikzpicture}[scale=0.9]
\begin{axis}[xlabel={Number of Appliances}, ylabel={Overall EA}, xmin=-1, xmax=18, ymin=-0.05, ymax=1.05, ytick={0, 0.20,0.40,0.60,0.80,1}, xtick={0, 2, 4, 6, 8, 10, 12, 14, 16, 18}, grid style=dashed, legend style={nodes={scale=0.55}, at={(0.25, 0.31)}, anchor=north east}]

\addplot[color=blue, mark=square] coordinates {(0,0.7681)(1,0.4917)(2,0.3891)(3,0.4927)(4,0.4709)(5,0.4772)(6,0.3905)(7,0.3599)(8,0.3274)(9,0.2947)(10,0.3022)(11,0.3173)(12,0.3766)(13,0.2961)(14,0.3345)(15,0.3351)(16,0.2831)(17,0.2788)};

\addplot[color=green, mark=o] coordinates {(0,0.7622)(1,0.6908)(2,0.6394)(3,0.6228)(4,0.6035)(5,0.5885)(6,0.5789)(7,0.5657)(8,0.5606)(9,0.5523)(10,0.5464)(11,0.5431)(12,0.5413)(13,0.5386)(14,0.5367)(15,0.5364)(16,0.5361)(17,0.536)};

\addplot[color=brown, mark=x] coordinates {(0,0.7131)(1,0.6386)(2,0.6309)(3,0.5973)(4,0.5979)(5,0.5914)(6,0.5743)(7,0.5744)(8,0.5704)(9,0.5704)(10,0.5701)(11,0.5693)(12,0.5694)(13,0.5668)(14,0.5643)(15,0.5644)(16,0.5649)(17,0.5646)};

\addplot[color=purple, mark=diamond] coordinates {(0,0.8165)(1,0.7373)(2,0.6821)(3,0.6559)(4,0.6357)(5,0.6185)(6,0.6109)(7,0.5975)(8,0.5911)(9,0.5812)(10,0.575)(11,0.5716)(12,0.5656)(13,0.5638)(14,0.5635)(15,0.5605)(16,0.5602)(17,0.5633)};

\addplot[color=red, mark=triangle] coordinates {(0,0.7474)(1,0.7521)(2,0.7160)(3,0.7128)(4,0.7104)(5,0.7094)(6,0.6888)(7,0.6772)(8,0.6652)(9,0.6642)(10,0.6536)(11,0.6536)(12,0.6446)(13,0.6441)(14,0.6435)(15,0.6440)(16,0.6344)(17,0.6345)};

\addplot[color=purple,mark=asterisk] coordinates {(0,0.9292)(1,0.9304)(2,0.9322)(3,0.9583)(4,0.9305)(5,0.9305)(6,0.9312)(7,0.9312)(8,0.9298)(9,0.9298)(10,0.9298)(11,0.9298)(12,0.93)(13,0.93)(14,0.93)(15,0.93)(16,0.93)(17,0.93)};

\legend{CO,SO,ALIP,STIP,SHMM,BQP}
\end{axis}
\end{tikzpicture}
    \caption{Overall Estimation Accuracy as the number of unmetered appliances varies.}
    \label{fig:synth1}
\end{figure}
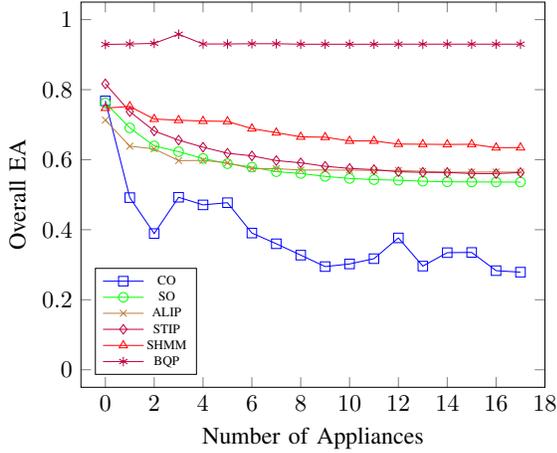

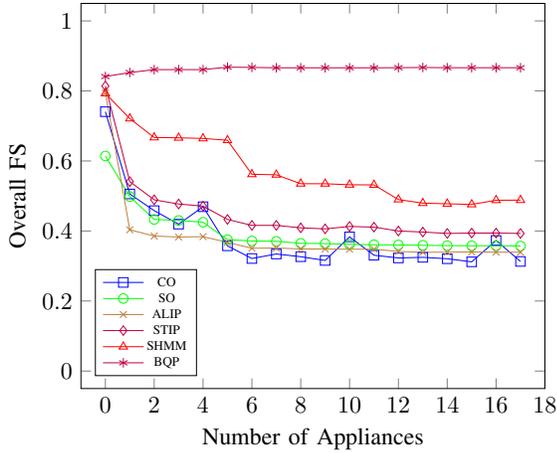
\begin{figure}
\centering
\begin{tikzpicture}[scale=0.9]
\begin{axis}[xlabel={Number of Appliances}, ylabel={Overall FS}, xmin=-1, xmax=18, ymin=-0.05, ymax=1.05, ytick={0, 0.20,0.40,0.60,0.80,1}, xtick={0, 2, 4, 6, 8, 10, 12, 14, 16, 18}, grid style=dashed, legend style={nodes={scale=0.55}, at={(0.25, 0.31)}, anchor=north east}]

\addplot[color=blue, mark=square] coordinates {(0,0.7407)(1,0.5057)(2,0.4581)(3,0.4195)(4,0.4693)(5,0.3576)(6,0.3214)(7,0.3343)(8,0.3265)(9,0.3158)(10,0.3831)(11,0.3307)(12,0.3230)(13,0.3247)(14,0.3207)(15,0.3119)(16,0.3729)(17,0.3130)};

\addplot[color=green, mark=o] coordinates {(0,0.6143)(1,0.4982)(2,0.4330)(3,0.4301)(4,0.4250)(5,0.3753)(6,0.3715)(7,0.3705)(8,0.3647)(9,0.3639)(10,0.3618)(11,0.3606)(12,0.3601)(13,0.3591)(14,0.3582)(15,0.3577)(16,0.3580)(17,0.3570)};

\addplot[color=brown, mark=x] coordinates {(0,0.7949)(1,0.4035)(2,0.3855)(3,0.3828)(4,0.3838)(5,0.3675)(6,0.3508)(7,0.3513)(8,0.3487)(9,0.3484)(10,0.3479)(11,0.3480)(12,0.3416)(13,0.3390)(14,0.3396)(15,0.3401)(16,0.3395)(17,0.3394)};

\addplot[color=purple,mark=diamond] coordinates {(0,0.814)(1,0.541)(2,0.489)(3,0.477)(4,0.471)(5,0.433)(6,0.416)(7,0.416)(8,0.409)(9,0.406)(10,0.413)(11,0.411)(12,0.400)(13,0.397)(14,0.393)(15,0.394)(16,0.394)(17,0.393)};

\addplot[color=red, mark=triangle] coordinates {(0,0.7930)(1,0.7211)(2,0.6670)(3,0.6661)(4,0.6643)(5,0.6591)(6,0.5619)(7,0.5605)(8,0.5350)(9,0.5349)(10,0.5320)(11,0.5316)(12,0.4898)(13,0.4795)(14,0.4777)(15,0.4758)(16,0.4878)(17,0.4881)};

\addplot[color=purple,mark=asterisk] coordinates {(0,0.8411)(1,0.8525)(2,0.8608)(3,0.8609)(4,0.8609)(5,0.8678)(6,0.8672)(7,0.8661)(8,0.8661)(9,0.8661)(10,0.8661)(11,0.8661)(12,0.8664)(13,0.8670)(14,0.8663)(15,0.8663)(16,0.8663)(17,0.8663)};

\legend{CO,SO,ALIP,STIP,SHMM,BQP}
\end{axis}
\end{tikzpicture}
    \caption{Overall FS as the number of unmetered appliances varies.}
    \label{fig:synth2}
\end{figure}

In Figures \ref{fig:synth1} and \ref{fig:synth2}, we report the average OEA and the average OFS obtained by all the methods when the number of appliances added to the aggregated signal increases. This picture shows the well-known difficulty encountered by most of the state-based approaches when the noise increases. Our method, on the other hand, does not suffer from the increase in noise, showing very high accuracy in all scenarios. This behavior stays essentially the same when looking at the single appliance. These results confirm the robustness to noise of our formulation, which allows overcoming the drawback of the other optimization-based approaches. This behavior is also confirmed by the experiments on real-world datasets described in the next section.

\subsection{Real-world Datasets}
We test our formulation on three real-world datasets: AMPDS \cite{makonin2013ampds}, UKDALE \cite{kelly2015uk} and REFIT \cite{murray2017electrical}, containing both the aggregate and the appliance-wise active power. From each dataset, we select at most two houses where the percentage of noise is up to 40\% and we use as target appliances the top-consuming ones. From AMPDS we select 6 appliances: clothes dryer (CD), dishwasher (DW), fan and thermostat (FT), entertainment (EN), fridge (FR), and heat pump (HP). We model FT and EN as always-on, whereas we model CD and HP as state machines. 
From UKDALE we select house 1 and 2. For house 1 we disaggregate 5 appliances: boiler (BO), washing machine (WM), dishwasher (DW), home theater (HT), and fridge (FR). For this house, WM and DW are treated as state machines, whereas BO and HT are treated as always-on. From House 2 we select 5 appliances: server (SE), kettle (KE), dishwasher (DW), fridge (FR), and modem (MO). Here, SE and MO are always-on, whereas DW is a state machine. 
Finally, for the REFIT dataset we choose house 3 and 9. For house 3 we use 4 appliances: fridge (FR), tumble dryer (TD), dishwasher (DW), and television (TV). We model TD and DW are state machines. For house 9 we select 4 appliances: fridge (FR), washing machine (WM), dishwasher (DW), and kettle (KT). We treat WM and DW are state machine appliances. In Table \ref{tab:real_dataset} we report the training, validation and test days and the percent-noisy measure for each house, computed according to \eqref{eq:noise}.

\begin{table}
\centering
\caption{Training, validation and test days.}
\scalebox{0.95}{
\begin{tabular}{|c|ccc|c|}
    \hline
    Dataset & Training & Validation & Test & Noise\\
    \hline
    \hline
    AMPDS & 10-23/04/12 & 24-30/04/12 & 01-07/05/12 & 38\%\\\hline
    UKDALE1 & 09-22/06/14 & 23-29/06/14 & 30-06/06-07/14 & 39\%\\
    UKDALE2 & 01-14/06/13 & 15-21/06/13 & 22-29/06/13 & 36\%\\\hline
    REFIT3 & 20-02/10-11/13 & 03-09/11/13 & 10-16/11/13 & 35\%\\
    REFIT9 & 18-01/11-12/14 & 02-08/12/14 & 09-15/12/14 & 40\%\\\hline
\end{tabular}}
\label{tab:real_dataset}
\end{table}

In Tables \ref{tab:ampd1}-\ref{tab:refit3}, we report for each method and for each house the Estimation Accuracy (EA) and  the Finite State F-Score (FS) for each appliance averaged on the test week. Furthermore, we report the average Overall EA and the average Overall FS. 
Note that we report the performance of our disaggregation algorithm with (BQPP) and without post-processing (BPQ), to evaluate the impact of the post-processing phase.

\begin{table*}
\footnotesize
\centering
\caption{AMPDS House 1}
 \begin{tabular}{|c|c|ccccccc|c|}
  \hline
    Method & Metric & Overall & CD & DW & FT & EN & FR & HP & Time\\
    \hline
    \hline
    \multirow{2}{*}{CO} & EA & 0.493 & 0.018 & -1.548 & 0.611 & 0.398 & 0.189 & 0.867 & \multirow{2}{*}{1133.5} \\
    & FS & 0.583 & 0.057 & 0.045 & 0.718 & 0.635 & 0.398 & 0.781 &\\
    \hline
    \multirow{2}{*}{SO} & EA & 0.569 & 0.815 & 0.500 & \textbf{0.978} & 0.195 & 0.490 & -0.141 & \multirow{2}{*}{1.8}\\
    & FS & 0.607 & 0.692 & 0.001 & \textbf{1.000} & 0.657 & 0.056 & 0.809 & \\ \hline
    \multirow{2}{*}{ALIP} & EA & 0.587 & 0.764 & 0.375 & $\textbf{0.978}$ & 0.361 & 0.195 & 0.629 & \multirow{2}{*}{3.40}\\
    & FS & 0.675 & 0.447 & 0.316 & \textbf{1.000} & 0.862 & 0.324 & 0.307 & \\
    \hline
    \multirow{2}{*}{STIP} & EA & 0.681 & 0.942 & 0.194 & \textbf{0.978} & 0.432 & 0.249 & 0.842 & \multirow{2}{*}{1206.6} \\
    & FS & 0.683 & 0.836 & 0.543 & $\textbf{1.000}$ & 0.872 & 0.513 & 0.653 & \\
    \hline
    \multirow{2}{*}{SHMM} & EA & 0.754 & 0.496 & 0.364 & 0.851 & 0.835 & 0.686 & 0.834 & \multirow{2}{*}{49.2}\\
    & FS & 0.786 & 0.310 & 0.561 & 0.841 & 0.943 & 0.662 & 0.757 & \\
    \hline
    \multirow{2}{*}{BQP} & EA & 0.889 & 0.955 & 0.863 & \textbf{0.978} & 0.849 & 0.721 & 0.907 & \multirow{2}{*}{52.5}\\
    & FS & \textbf{0.942} & \textbf{0.932} & \textbf{0.844} & \textbf{1.000} & \textbf{0.978} & \textbf{0.760} & \textbf{0.919} & \\\hline
    \multirow{2}{*}{BQPP} & EA & \textbf{0.902} & \textbf{0.976} & \textbf{0.888} & \textbf{0.978} & \textbf{0.862} & \textbf{0.736} & \textbf{0.923} & \multirow{2}{*}{52.8}\\
    & FS & \textbf{0.942} & \textbf{0.932} & \textbf{0.844} & \textbf{1.000} & \textbf{0.978} & \textbf{0.760} & \textbf{0.919} & \\\hline
\end{tabular}
\label{tab:ampd1}
\end{table*}

\begin{table*}
\footnotesize
\centering
\caption{UKDALE Houses 1 and 2} 
\begin{tabular}{|c|c|cccccc|c||cccccc|c|}
    \hline
    & & \multicolumn{7}{|c|}{House 1} &     \multicolumn{7}{ |c|}{House 2} \\
    \hline
    Method & Metric & Overall & BO & WM & DW & HT & FR & Time & Overall & SE & KE & DW & FR & MO & Time \\
    \hline
    \hline
    \multirow{2}{*}{CO} & EA & -0.192 & 0.375 & -0.501 & -2.029 & -1.421 & 0.411 & \multirow{2}{*}{1034.5} & 0.403 & 0.754 & 0.622 & -0.030 & 0.584 & 0.742 & \multirow{2}{*}{755.7}\\
    & FS & 0.421 & 0.587 & 0.155 & 0.074 & 0.598 & 0.422 & & 0.662 & 0.710 & 0.535 & 0.101 & 0.786 & 0.755 &\\
    \hline
    \multirow{2}{*}{SO} & EA & 0.305 & 0.977 & 0.485 & -3.221 & 0.559 & 0.574 & \multirow{2}{*}{1.0} & 0.640 & 0.990 & 0.507 & 0.589 & 0.620 & \textbf{0.997} & \multirow{2}{*}{0.9}\\
    & FS & 0.677 & 0.997 & 0.480 & 0.105 & 0.829 & 0.677 & & 0.406 & \textbf{1.000} & 0.009 & 0.003 & 0.012 & \textbf{1.000} &\\
    \hline
    \multirow{2}{*}{ALIP} & EA & 0.174 & 0.195 & 0.383 & -1.468 & -0.933 & 0.557 & \multirow{2}{*}{9.5} & 0.564 & 0.990 & 0.661 & 0.325 & 0.659 & \textbf{0.997} & \multirow{2}{*}{3.2}\\
    & FS & 0.561 & 0.789 & 0.442 & 0.149 & 0.616 & 0.454 &  & 0.638 & \textbf{1.000} & 0.105 & 0.209 & 0.597 & \textbf{1.000} & \\
    \hline
    \multirow{2}{*}{STIP} & EA & 0.157 & 0.313 & 0.581 & -0.506 & -0.651 & 0.569 & \multirow{2}{*}{723.6} & 0.579 & 0.990 & 0.661 & 0.387 & 0.496 & \textbf{0.997} & \multirow{2}{*}{957.1}\\
    & FS & 0.579 & 0.821 & 0.405 & 0.178 & 0.843 & 0.549 & & 0.778 & \textbf{1.000} & 0.565 & 0.236 & 0.857 & \textbf{1.000} &\\
    \hline
    \multirow{2}{*}{SHMM} & EA & 0.629 & 0.884 & 0.478 & 0.346 & 0.184 & 0.731 & \multirow{2}{*}{14.7} & 0.723 & 0.966 & 0.519 & 0.571 & 0.802 & 0.981 & \multirow{2}{*}{4.8}\\
    & FS & 0.782 & 0.913 & 0.306 & 0.779 & 0.818 & 0.712 & & 0.856 & 0.982 & 0.780 & 0.382 & 0.876 & 0.982 & \\ \hline
    \multirow{2}{*}{BQP} & EA & 0.869 & 0.958 & 0.711 & 0.903 & 0.699 & 0.795 & \multirow{2}{*}{25.4} & 0.909 & 0.990 & 0.778 & 0.907 & 0.921 & \textbf{0.997} & \multirow{2}{*}{5.9}\\
    & FS & \textbf{0.951} & \textbf{0.997} & \textbf{0.884} & \textbf{0.802} & \textbf{0.967} & \textbf{0.853} & & \textbf{0.970} & \textbf{1.000} & \textbf{0.789} & \textbf{0.868} & \textbf{0.975} & \textbf{1.000} & \\
    \hline
    \multirow{2}{*}{BQPP} & EA & \textbf{0.894} & \textbf{0.979} & \textbf{0.746} & \textbf{0.932} & \textbf{0.702} & \textbf{0.809} & \multirow{2}{*}{25.7} & \textbf{0.926} & \textbf{0.991} & \textbf{0.790} & \textbf{0.933} & \textbf{0.938} & \textbf{0.997} & \multirow{2}{*}{6.1}\\
    & FS & \textbf{0.951} & \textbf{0.997} & \textbf{0.884} & \textbf{0.802} & \textbf{0.967} & \textbf{0.853} & & \textbf{0.970} & \textbf{1.000} & \textbf{0.789} & \textbf{0.868} & \textbf{0.975} & \textbf{1.000} & \\
    \hline
\end{tabular}
\label{tab:uk1}
\end{table*}

\begin{table*}[!ht]
\footnotesize
\centering
\caption{REFIT Houses 3 and 9}
\begin{tabular}{|c|c|ccccc|c||ccccc|c|}
    \hline
    & & \multicolumn{6}{|c|}{House 3} &     \multicolumn{6}{ |c|}{House 9} \\
    \hline
 Method & Metric & Overall & FR & TD & DW & TV & Time & Overall & FR & WM & DW & KE & Time\\
    \hline
    \hline
    \multirow{2}{*}{CO} & EA & 0.264 & 0.482 & 0.344 & -0.092 & 0.589 & \multirow{2}{*}{989.1} & -0.262 & 0.358 & -0.467 & -0.017 & -0.043 & \multirow{2}{*}{687.3}\\
    & FS & 0.501 & 0.676 & 0.441 & 0.160 & 0.426 & & 0.164 & 0.240 &  0.074 &  0.116 &  0.000 &\\
    \hline
    \multirow{2}{*}{SO} & EA & 0.375 & 0.640 & 0.281 & -0.863 & 0.763 & \multirow{2}{*}{0.9} & 0.166 & 0.619 & 0.522 & 0.471 & -1.007 & \multirow{2}{*}{1.1}\\
    & FS & 0.630 & 0.578 & 0.409 & 0.218 & 0.805 & & 0.603 & 0.702 &  0.106 &  0.282 &  0.332 &\\
    \hline
    \multirow{2}{*}{ALIP} & EA & 0.538 & 0.535 & 0.554 & 0.278 & 0.874 &  \multirow{2}{*}{1.1} & 0.226 & 0.658 & 0.402 & 0.432 & 0.486 & \multirow{2}{*}{0.9}\\
    & FS & 0.673 & 0.548 & 0.465 & 0.350 & 0.872 & & 0.490  & 0.587 &  0.120 &  0.204  & 0.000 &\\
    \hline
    \multirow{2}{*}{STIP} & EA & 0.553 & 0.495 & 0.669 & 0.131 & 0.752 & \multirow{2}{*}{468.9} & 0.556 & 0.667 & 0.398 & 0.469 & 0.581 & \multirow{2}{*}{509.6}\\
    & FS & 0.526 & 0.549 & 0.637 & 0.124 & 0.740 & & 0.525 & 0.599 & 0.197 & 0.308 & 0.416 &\\
    \hline
    \multirow{2}{*}{SHMM} & EA & 0.676 & 0.798 & 0.494 & 0.500 & 0.889 & \multirow{2}{*}{3.1} & 0.695 & 0.762 & 0.485 & 0.745 & 0.500 & \multirow{2}{*}{4.9}\\
    & FS & 0.744 & 0.802 & 0.103 & 0.358 & 0.782 & & 0.637 & 0.789 & 0.272 & 0.611 & 0.005 &\\
    \hline
    \multirow{2}{*}{BQP} & EA & 0.828 & 0.807 & 0.841 & 0.819 & 0.881 & \multirow{2}{*}{3.6} & 0.806 & 0.771 & 0.737 & 0.835 & 0.762 & \multirow{2}{*}{7.8}\\
    & FS & \textbf{0.862} & \textbf{0.858} & \textbf{0.855} & \textbf{0.815} & \textbf{0.892} & & \textbf{0.801} & \textbf{0.814} & \textbf{0.771} & \textbf{0.799} & \textbf{0.643} & \\
    \hline
    \multirow{2}{*}{BQPP} & EA & \textbf{0.847} & \textbf{0.813} & \textbf{0.878} & \textbf{0.844} & \textbf{0.894} & \multirow{2}{*}{3.8} & \textbf{0.823} & \textbf{0.788} & \textbf{0.744} & \textbf{0.846} & \textbf{0.772} & \multirow{2}{*}{7.9}\\
    & FS & \textbf{0.862} & \textbf{0.858} & \textbf{0.855} & \textbf{0.815} & \textbf{0.892} & & \textbf{0.801} & \textbf{0.814} & \textbf{0.771} & \textbf{0.799} & \textbf{0.643} & \\
    \hline
\end{tabular}
\label{tab:refit3}
\end{table*}
Finally, we report the average execution time in seconds for each method. 
All the compared approaches are based on a two-stage procedure: an offline phase for parameter estimation and an online phase for the actual disaggregation. In all the cases, the reported time represents the computational time of the online procedure only. Note that, anyway, the offline phase of our approach is extremely cheap in terms of computational time: the most expensive part is the grid search for $\lambda_1$ and $\lambda_2$ that requires few minutes. Looking at the execution times, the largest ones are required by CO, due to the larger feasible region, and by STIP, due to the nonlinear constraints induced by the state machine constraints. The solution of our BQP problem only requires a few seconds, despite the large size: indeed, we have thousands of binary variables and thousands of constraints.
Therefore, our approach could be used in a real-time setting, at least for the considered data granularity, and the considered scheduling horizon. This makes the proposed method suitable for low-cost microcontrollers.

\begin{figure}[!ht]
\centering
\begin{tikzpicture}[scale=0.9]
\begin{axis}[title={}, xlabel={Houses}, ylabel={Overall FS}, xmin=-1, xmax=9, ymin=0.05, ymax=1.10, ticklabel style={font=\footnotesize}, ytick={0.0, 0.20,0.40,0.60,0.80,1},
xtick={0, 2, 4, 6, 8}, xticklabels={AMPDS, UKDALE1, UKDALE2, REFIT3, REFIT9}, grid style=dashed, legend style={nodes={scale=0.50}, at={(0.24, 0.29)}, anchor=north east}]

\addplot[color=blue, mark=square] coordinates
{(0,0.583)(2,0.421)(4,0.662)(6,0.501)(8,0.164)};

\addplot[color=green, mark=o] coordinates 
{(0,0.607)(2,0.677)(4,0.406)(6,0.63)(8,0.603)};

\addplot[color=brown, mark=x] coordinates 
{(0,0.675)(2,0.561)(4,0.638)(6,0.673)(8,0.49)};

\addplot[color=yellow, mark=diamond] coordinates 
{(0,0.683)(2,0.579)(4,0.778)(6,0.526)(8,0.525)};

\addplot[color=red, mark=triangle] coordinates
{(0,0.786)(2,0.782)(4,0.856)(6,0.744)(8,0.637)};

\addplot[color=purple,mark=asterisk] coordinates 
{(0,0.942)(2,0.951)(4,0.970)(6,0.862)(8,0.801)};

\legend{CO,SO,ALIP,STIP,SHMM,BQP}
\end{axis}
\end{tikzpicture}
    \caption{Overall FS-Fscore for each house and method.}
    \label{fig:fs_all}
\end{figure}
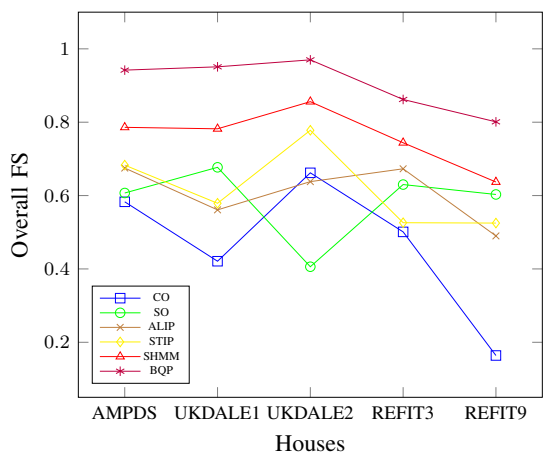

\begin{figure*}
  \centering
  \includegraphics[width=\textwidth]{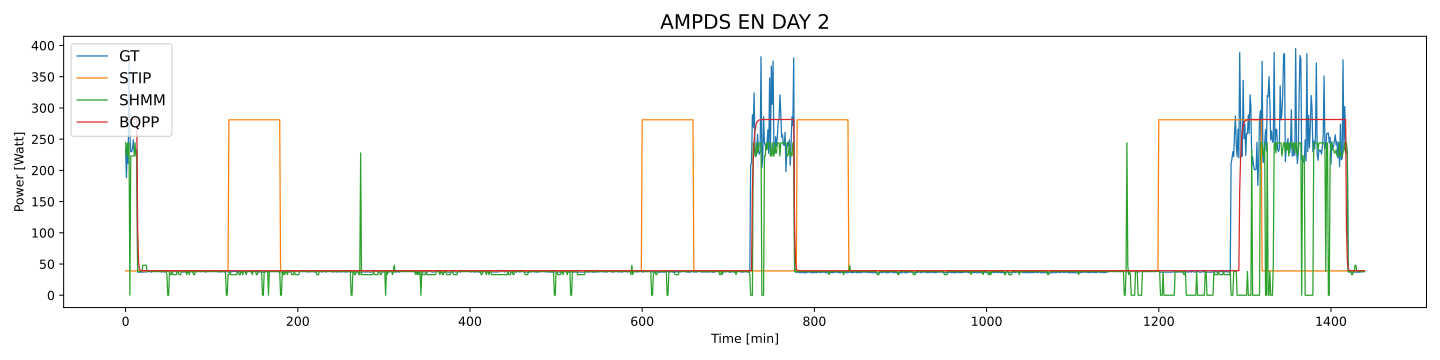}
  \includegraphics[width=\textwidth]{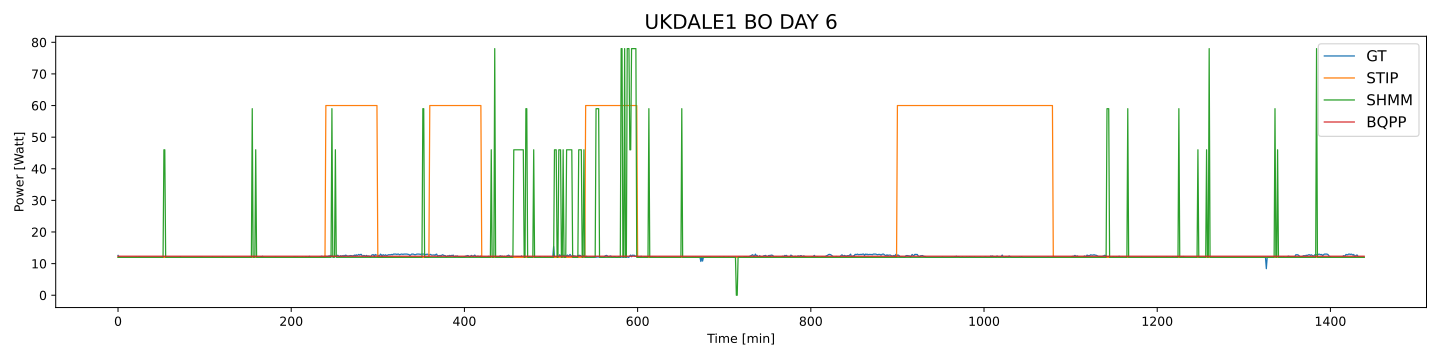}
  \includegraphics[width=\textwidth]{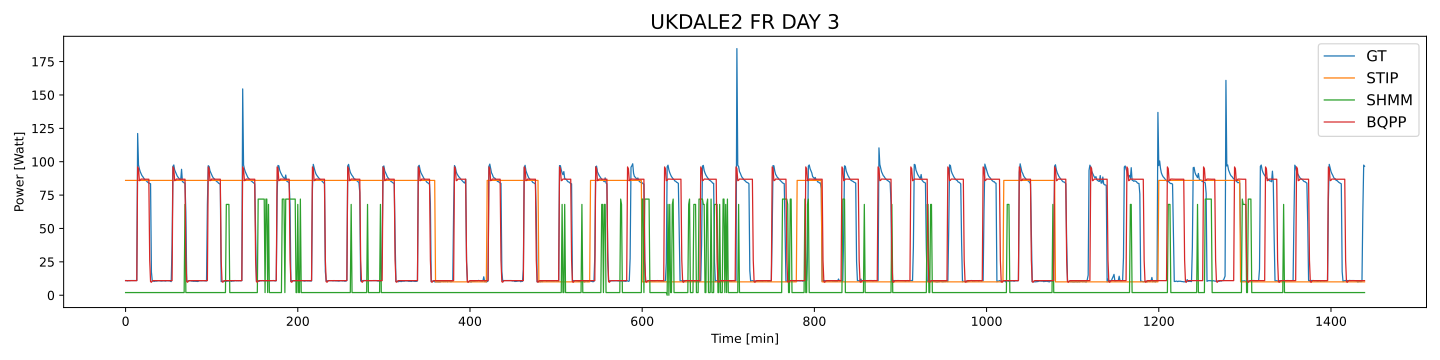}
  \includegraphics[width=\textwidth]{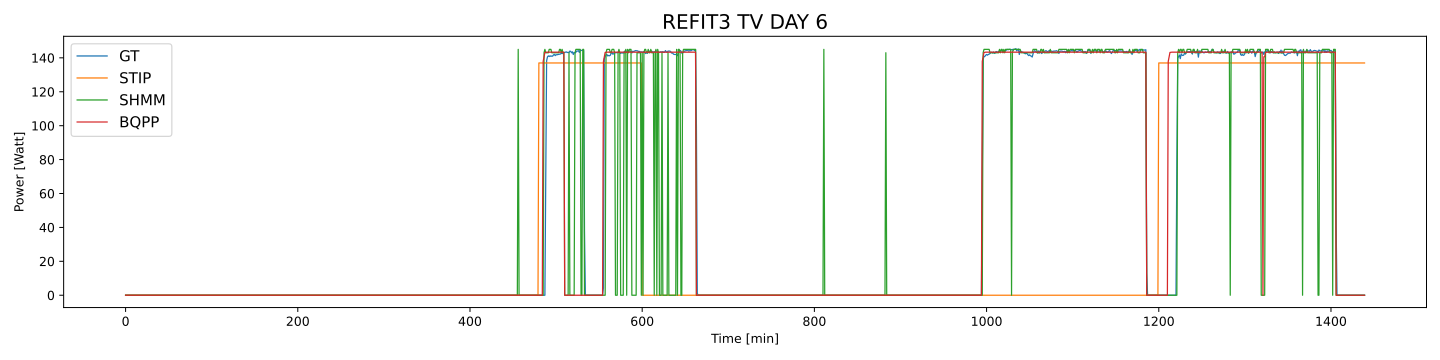}
  \includegraphics[width=\textwidth]{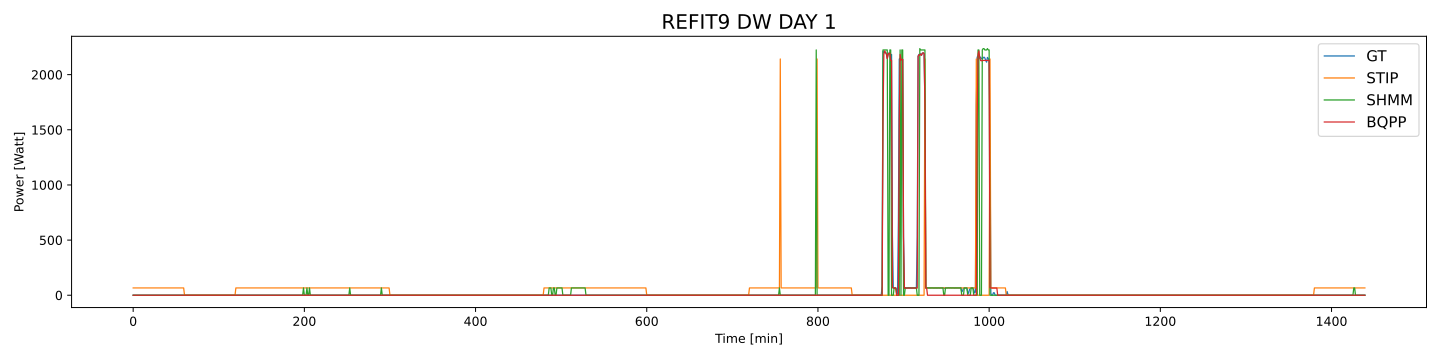}
  \caption{Ground truth (GT) and estimated power consumption of different appliances.}
  \label{fig:one_day_disaggregation}
\end{figure*}

As for the performance metrics, it turns out that our method is always the best in both metrics on all the appliances. To get a better picture of the different methods' behavior, we extract from the tables the average Overall EA and the average Overall FS on all the different houses, and plot them in Figures \ref{fig:ea_all} and \ref{fig:fs_all}. The superiority of our method in both metrics is evident, and in Figure \ref{fig:ea_all} also emerges the advantage of the post-processing phase, which has an impact only on the EA since it does not influence the activations. Looking at the detailed performance on the single appliances in Tables \ref{tab:ampd1}-\ref{tab:refit3}, it can be seen that most of the approaches are able to correctly disaggregate always-on appliances (see for example appliance FT of AMPDS), whereas the BQP is by far more robust on appliances that have a variable number of activations (see for example appliances DW of both house 1 and 2 of UKDALE). The high performance on these appliances derives from the ability of the proposed formulation to disambiguate the appliances' states, combined with the correct estimation of the power states. The beneficial effect of the post-processing is confirmed since the EA metric improves on all the appliances in all the datasets. 
To better show the performance of the proposed algorithm, in Figure \ref{fig:one_day_disaggregation} we show the 1-day disaggregation results produced by the top three methods BQPP, SHMM and STIP and we select appliances with different power signatures and behaviours. We can see that the fridge typically shows a peak at the beginning of its working cycle. Differently from STIP that produces a flat estimate, our approach manages to capture this behaviour thanks to the AR coefficients. On the other hand, the entertainment is a Type-III appliance whose activations are successfully captured even though the consumption is approximated by the average power. Looking at the dishwasher, we can see that our method is successful at capturing with very good accuracy the activations and the kW-level power consumption.
In general, our approach is able to correctly disambiguate appliances with similar consumption patterns and to correctly estimate the load of active appliances in a real-world scenario. The competitors tend to produce false positives for the appliance that are not active (see the disaggregation results for television and boiler in Figure \ref{fig:one_day_disaggregation}). Furthermore, the combination of our post-processing phase and of our parameter estimation procedure allows to capture with great accuracy the scale of the active power states, increasing the EA metric. 

\begin{figure}[!ht]
\centering
\begin{tikzpicture}[scale=0.9]
\begin{axis}[title={}, xlabel={Houses}, ylabel={Overall EA}, xmin=-1, xmax=9, ymin=-0.60, ymax=1.10, ticklabel style={font=\footnotesize}, ytick={-0.40, -0.20,0.0, 0.20,0.40,0.60,0.80,1,1.20},
xtick={0, 2, 4, 6, 8}, xticklabels={AMPDS, UKDALE1, UKDALE2, REFIT3, REFIT9}, grid style=dashed, legend style={nodes={scale=0.50}, at={(0.24, 0.33)}, anchor=north east}]

\addplot[color=blue, mark=square] coordinates
{(0,0.493)(2,-0.192)(4,0.403)(6,0.264)(8,-0.262)};

\addplot[color=green, mark=o] coordinates 
{(0,0.569)(2,0.305)(4,0.64)(6,0.375)(8,0.166)};

\addplot[color=brown, mark=x] coordinates 
{(0,0.587)(2,0.174)(4,0.564)(6,0.538)(8,0.226)};

\addplot[color=yellow, mark=diamond] coordinates 
{(0,0.681)(2,0.157)(4,0.579)(6,0.553)(8,0.556)};

\addplot[color=red, mark=triangle] coordinates
{(0,0.754)(2,0.629)(4,0.723)(6,0.676)(8,0.695)};

\addplot[color=purple,mark=asterisk] coordinates 
{(0,0.889)(2,0.869)(4,0.909)(6,0.828)(8,0.806)};

\addplot[color=purple,mark=diamond] coordinates 
{(0,0.902)(2,0.894)(4,0.926)(6,0.847)(8,0.823)};

\legend{CO,SO,ALIP,STIP,SHMM,BQP,BQPP}
\end{axis}
\end{tikzpicture}
    \caption{Overall Estimation Accuracy for each house and method.}
    \label{fig:ea_all}
\end{figure}
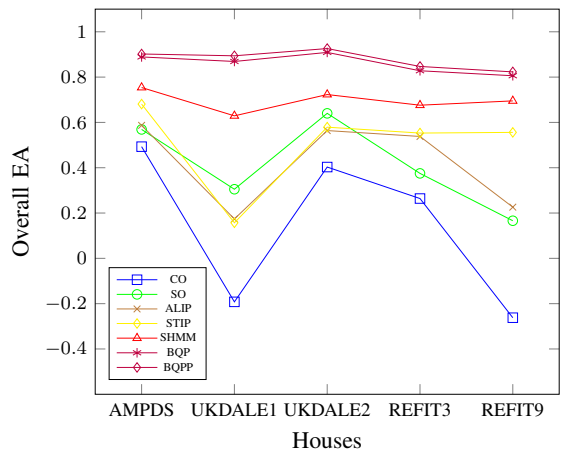

\subsection{Parameters Stability}
\label{sec:par_upd}
The proposed optimization model relies on several parameters that may need to be updated in order to reflect changes in the system. As long as there is no significant change in the power profiles, the estimated parameters should remain stable over time.
In the following, we conduct experiments to study whether the freshness and staleness of the parameters impact the accuracy of the disaggregation. We compare two different scenarios:
\paragraph{Scenario A (Fixed)} the parameters are estimated just once, at the beginning of the experiment, using the training and validation weeks reported in Table \ref{tab:real_dataset}. Then,
we simulate a periodic disaggregation once per month until six months after the beginning of the experiment. Therefore,
we perform Stage 2 of our disaggregation procedure on six different test weeks, further and further away from the weeks used for parameter estimation. In this setting, we evaluate what happens when the parameters are not updated for a long time.
\paragraph{Scenario B (Updated)} the overall disaggregation procedure is repeated at one month interval up to six months ahead on the same test weeks of Scenario A. However, in this case the parameters are re-estimated by using two weeks for training and one for validation immediately before the test week. This is an ``ideal scenario'' where Stage 1 is executed right before Stage 2. We note that this could still be viable in practice if historical observations were available.

\begin{figure*}
\begin{tikzpicture}
    \begin{axis}[
        legend style={nodes={scale=0.45}},
        legend columns=2, 
        xtick={1, 2, 3, 4, 5, 6},
        xticklabels={},
        ytick={0.80, 0.85, 0.90},
        title={CD},
        ylabel={EA},
        width=0.34\textwidth,
        height=3.25cm,
        legend pos=north west,
        ymin = 0.75,
        ymax = 0.95
    ]
        \addplot [blue, mark=o]  coordinates {(1,0.846)(2,0.871)(3,0.816)(4,0.855)(5,0.878)(6,0.841)};
        \addplot [red, mark=x]  coordinates {(1,0.842)(2,0.87)(3,0.816)(4,0.859)(5,0.878)(6,0.841)};
        \legend{Fixed, Updated}
    \end{axis}
\end{tikzpicture}
\begin{tikzpicture}
   \begin{axis}[
        legend style={nodes={scale=0.45}},
        legend columns=2, 
        xtick={1, 2, 3, 4, 5, 6},
        xticklabels={},
        ytick={0.75, 0.80, 0.85},
        title={DW},
        ylabel={EA},
        width=0.34\textwidth,
        height=3.25cm,
        legend pos=north west,
        ymin = 0.70,
        ymax = 0.90
    ]
        \addplot [blue, mark=o]  coordinates {(1,0.822)(2,0.746)(3,0.774)(4,0.828)(5,0.792)(6,0.764)};
        \addplot [red, mark=x]  coordinates {(1,0.832)(2,0.739)(3,0.777)(4,0.828)(5,0.795)(6,0.766)};
        \legend{Fixed, Updated}
    \end{axis}
\end{tikzpicture}
\begin{tikzpicture}
    \begin{axis}[
        legend style={nodes={scale=0.45}},
        legend columns=2, 
        xtick={1, 2, 3, 4, 5, 6},
        xticklabels={},
        ytick={0.90, 0.95, 1.0},
        title={FT},
        ylabel={EA},
        width=0.34\textwidth,
        height=3.25cm,
        legend pos=north west,
        ymin = 0.85,
        ymax = 1.05
    ]
        \addplot [blue, mark=o]  coordinates {(1,0.947)(2,0.961)(3,0.981)(4,0.898)(5,0.958)(6,0.906)};
        \addplot [red, mark=x]  coordinates {(1,0.94)(2,0.966)(3,0.981)(4,0.898)(5,0.968)(6,0.923)};
        \legend{Fixed, Updated}
    \end{axis}
\end{tikzpicture}
\begin{tikzpicture}
    \begin{axis}[
        legend style={nodes={scale=0.45}},
        legend columns=2, 
        xtick={1, 2, 3, 4, 5, 6},
        xticklabels={1, 2, 3, 4, 5, 6},
        xlabel={Months Ahead},
        ytick={0.80, 0.85, 0.90},
        title={},
        ylabel={FS},
        width=0.34\textwidth,
        height=3.25cm,
        legend pos=north west,
        ymin = 0.75,
        ymax = 0.95
    ]
        \addplot [blue, mark=o]  coordinates {(1,0.832)(2,0.828)(3,0.85)(4,0.86)(5,0.824)(6,0.865)};
        \addplot [red, mark=x]  coordinates {(1,0.84)(2,0.832)(3,0.846)(4,0.862)(5,0.82)(6,0.863)};
        \legend{Fixed, Updated}
    \end{axis}
\end{tikzpicture}
\begin{tikzpicture}
     \begin{axis}[
        legend style={nodes={scale=0.45}},
        legend columns=2, 
        xtick={1, 2, 3, 4, 5, 6},
        xticklabels={1, 2, 3, 4, 5, 6},
        xlabel={Months Ahead},
        ytick={0.75, 0.8, 0.85},
        title={},
        ylabel={FS},
        width=0.34\textwidth,
        height=3.25cm,
        legend pos=south west,
        ymin = 0.70,
        ymax = 0.90
    ]
        \addplot [blue, mark=o]  coordinates {(1,0.838)(2,0.825)(3,0.785)(4,0.839)(5,0.799)(6,0.794)};
        \addplot [red, mark=x]  coordinates {(1,0.836)(2,0.826)(3,0.785)(4,0.844)(5,0.798)(6,0.799)};
        \legend{Fixed, Updated}
    \end{axis}
\end{tikzpicture}
\begin{tikzpicture}
    \begin{axis}[
        legend style={nodes={scale=0.45}},
        legend columns=2, 
        xtick={1, 2, 3, 4, 5, 6},
        xticklabels={1, 2, 3, 4, 5, 6},
        xlabel={Months Ahead},
        ytick={0.90, 0.95, 1.00},
        title={},
        ylabel={FS},
        width=0.34\textwidth,
        height=3.25cm,
        legend pos=south west,
        ymin = 0.85,
        ymax = 1.05,
    ]
        \addplot [blue, mark=o]  coordinates {(1,1.0)(2,1.0)(3,0.999)(4,1.0)(5,0.988)(6,0.983)};
        \addplot [red, mark=x]  coordinates {(1,1.0)(2,1.0)(3,1.0)(4,1.0)(5,1.0)(6,0.99)};
        \legend{Fixed, Updated}
    \end{axis}
\end{tikzpicture}
\begin{tikzpicture}
    \begin{axis}[
        legend style={nodes={scale=0.45}},
        legend columns=2, 
        xtick={1, 2, 3, 4, 5, 6},
        xticklabels={},
        ytick={0.50, 0.75, 1.0},
        title={EN},
        ylabel={EA},
        width=0.34\textwidth,
        height=3.25cm,
        legend pos=south west,
        ymin = 0.30,
        ymax = 1.10
    ]
        \addplot [blue, mark=o]  coordinates {(1,0.910)(2,0.811)(3,0.757)(4,0.659)(5,0.573)(6,0.623)};
        \addplot [red, mark=x]  coordinates {(1,0.915)(2,0.821)(3,0.837)(4,0.829)(5,0.843)(6,0.869)};
        \legend{Fixed, Updated}
    \end{axis}
\end{tikzpicture}
\begin{tikzpicture}
   \begin{axis}[
        legend style={nodes={scale=0.45}},
        legend columns=2, 
        xtick={1, 2, 3, 4, 5, 6},
        xticklabels={},
        ytick={0.70, 0.75, 0.80},
        title={FR},
        ylabel={EA},
        width=0.34\textwidth,
        height=3.25cm,
        legend pos=north west,
        ymin = 0.65,
        ymax = 0.875
    ]
        \addplot [blue, mark=o]  coordinates {(1,0.737)(2,0.72)(3,0.717)(4,0.804)(5,0.719)(6,0.723)};
        \addplot [red, mark=x]  coordinates {(1,0.742)(2,0.717)(3,0.713)(4,0.8)(5,0.723)(6,0.724)};
        \legend{Fixed, Updated}
    \end{axis}
\end{tikzpicture}
\begin{tikzpicture}
    \begin{axis}[
        legend style={nodes={scale=0.45}},
        legend columns=2, 
        xtick={1, 2, 3, 4, 5, 6},
        xticklabels={},
        ytick={0.75, 0.85, 0.95},
        title={HP},
        ylabel={EA},
        width=0.34\textwidth,
        height=3.25cm,
        legend pos=south east,
        ymin = 0.65,
        ymax = 1.05
    ]
        \addplot [blue, mark=o]  coordinates { (1,0.821)(2,0.834)(3,0.897)(4,0.963)(5,0.904)(6,0.817)};
        \addplot [red, mark=x]  coordinates { (1,0.82)(2,0.834)(3,0.899)(4,0.964)(5,0.932)(6,0.879)};
        \legend{Fixed, Updated}
    \end{axis}
\end{tikzpicture}
\begin{tikzpicture}
    \begin{axis}[
        legend style={nodes={scale=0.45}},
        legend columns=2, 
        xtick={1, 2, 3, 4, 5, 6},
        xticklabels={1, 2, 3, 4, 5, 6},
        xlabel={Months Ahead},
        ytick={0.50, 0.75, 1.0},
        title={},
        ylabel={FS},
        width=0.34\textwidth,
        height=3.25cm,
        legend pos=south west,
        ymin = 0.40,
        ymax = 1.10
    ]
        \addplot [blue, mark=o]  coordinates {(1,0.99)(2,0.967)(3,0.897)(4,0.786)(5,0.642)(6,0.647)};
        \addplot [red, mark=x]  coordinates {(1,1.0)(2,0.967)(3,0.967)(4,0.936)(5,0.922)(6,0.917)};
        \legend{Fixed, Updated}
    \end{axis}
\end{tikzpicture}
\begin{tikzpicture}
     \begin{axis}[
        legend style={nodes={scale=0.45}},
        legend columns=2, 
        xtick={1, 2, 3, 4, 5, 6},
        xticklabels={1, 2, 3, 4, 5, 6},
        xlabel={Months Ahead},
        ytick={0.75, 0.80, 0.85},
        title={},
        ylabel={FS},
        width=0.34\textwidth,
        height=3.25cm,
        legend pos=north west,
        ymin = 0.65,
        ymax = 0.95
    ]
        \addplot [blue, mark=o]  coordinates {(1,0.786)(2,0.725)(3,0.713)(4,0.834)(5,0.756)(6,0.775)};
        \addplot [red, mark=x]  coordinates {(1,0.782)(2,0.718)(3,0.719)(4,0.835)(5,0.751)(6,0.782)};
        \legend{Fixed, Updated}
    \end{axis}
\end{tikzpicture}
\begin{tikzpicture}
     \begin{axis}[
        legend style={nodes={scale=0.45}},
        legend columns=2, 
        xtick={1, 2, 3, 4, 5, 6},
        xticklabels={1, 2, 3, 4, 5, 6},
        xlabel={Months Ahead},
        ytick={0.75, 0.85, 0.95},
        title={},
        ylabel={FS},
        width=0.34\textwidth,
        height=3.25cm,
        legend pos=south east,
        ymin = 0.65,
        ymax = 1.05
    ]
        \addplot [blue, mark=o]  coordinates {(1,0.798)(2,0.839)(3,0.916)(4,0.977)(5,0.934)(6,0.838)};
        \addplot [red, mark=x]  coordinates {(1,0.801)(2,0.833)(3,0.91)(4,0.985)(5,0.977)(6,0.882)};
        \legend{Fixed, Updated}
    \end{axis}
\end{tikzpicture}
\caption{Appliance-wise accuracy of the proposed approach in Scenario A (Fixed Parameters) and Scenario B (Updated Parameters).}
\label{fig:stability_app}
\end{figure*}
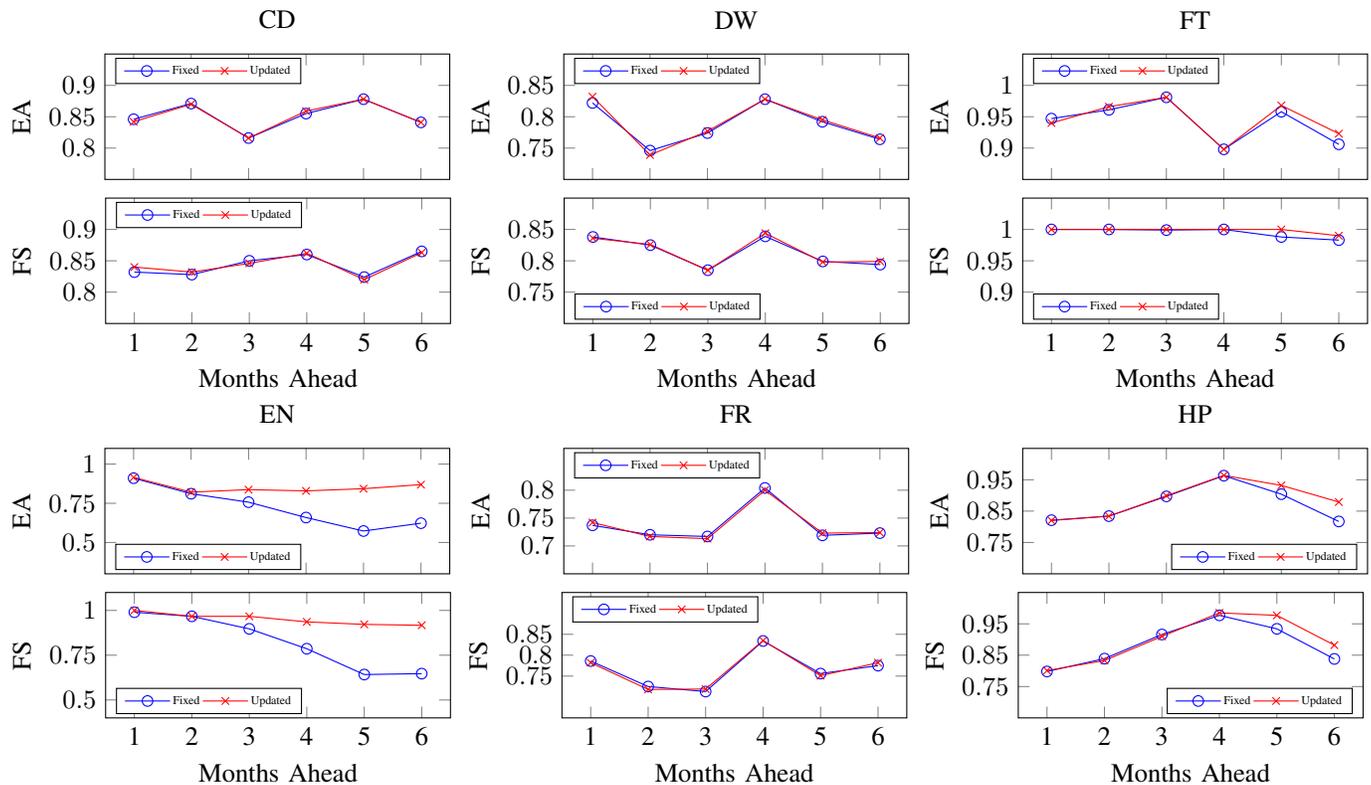

To conduct the experiments, we use the AMPDS dataset since it has a larger number of appliances, and thus a larger number of parameters. Results are reported in Figures \ref{fig:stability_overall} and \ref{fig:stability_app}. In Figure \ref{fig:stability_overall} we evaluate in both scenarios the overall EA and FS scores obtained in the six test weeks placed at one month interval. Interestingly enough the overall performance do not differ too much in the two settings, with a slight decrease in Scenario A starting from the fourth month. To better understand the two performance, in Figure \ref{fig:stability_app} the same metrics are reported in both scenarios on the single appliances. We can see that CD, DW, FT, and FR are not affected by the parameter staleness, since the accuracy in the two scenarios is almost identical. On the other hand, the only two appliances that suffer from staleness are the HP starting from the fifth month and EN starting form the third one. According to the data sheet in \cite{makonin2013ampds}, the heat pump cools the house in summer and heats the house in winter. The behavioural change of HP is related to a seasonal change, and hence the need to use the appliance more often in the heating setting. Looking at the ground truth signature, this change implies the presence of an additional state that is successfully captured in Scenario B when the parameters are refreshed, but not in Scenario A.
As for the EN, in \cite{makonin2013ampds} is treated as a composite appliance having multiple loads such as television set, personal video recorder and external amplifier. Looking at the active power signature it seems that from the fourth month onwards, it continuously changes its behaviour displaying a much bigger variability that is better captured in Scenario B and hence justifying the need for a new parameter estimation.

\begin{figure}[!ht]
    \centering
    \begin{tikzpicture}
    \begin{axis}[
        legend style={nodes={scale=0.45}},
        legend columns=2, 
        xtick={1, 2, 3, 4, 5, 6},
        xticklabels={},
        ytick={0.80, 0.85, 0.90},
        title={Overall},
        ylabel={EA},
        width=0.34\textwidth,
        height=3.25cm,
        legend pos=south east,
        ymin = 0.75,
        ymax = 0.95
    ]
        \addplot [blue, mark=o]  coordinates { (1,0.845)(2,0.817)(3,0.816)(4,0.872)(5,0.873)(6,0.844)};
        \addplot [red, mark=x]  coordinates { (1,0.843)(2,0.821)(3, 0.828)(4, 0.889)(5, 0.897)(6, 0.876)};
        \legend{Fixed, Updated}
    \end{axis}
    \end{tikzpicture}
    \begin{tikzpicture}
    \begin{axis}[
        legend style={nodes={scale=0.45}},
        legend columns=2, 
        xtick={1, 2, 3, 4, 5, 6},
        xticklabels={1, 2, 3, 4, 5, 6},
        xlabel={Months Ahead},
        ytick={0.80, 0.85, 0.90},
        title={},
        ylabel={FS},
        width=0.34\textwidth,
        height=3.25cm,
        legend pos=south west,
        ymin = 0.75,
        ymax = 0.95
    ]
        \addplot [blue, mark=o]  coordinates {(1,0.898) (2,0.887)(3, 0.846)(4, 0.869)(5, 0.846)(6, 0.831)};
        \addplot [red, mark=x]  coordinates {(1,0.900)(2, 0.889)(3, 0.856)(4, 0.882)(5, 0.876)(6, 0.868)};
        \legend{Fixed, Updated}
    \end{axis}
\end{tikzpicture}
    \caption{Overall accuracy of the proposed approach in Scenario A (Fixed Parameters) and Scenario B (Updated Parameters).}
    \label{fig:stability_overall}
\end{figure}
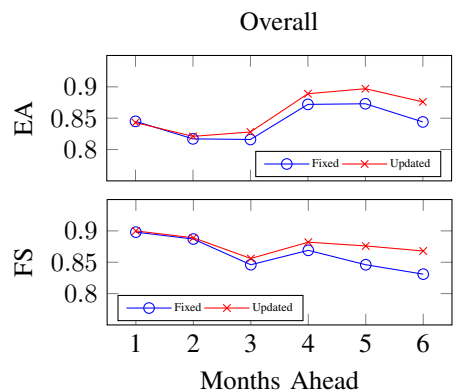

\section{Conclusions}
In this paper, we have proposed a novel optimization-based algorithm for NILM. Our approach can separate many appliances almost perfectly using only the total aggregate signal. It is computationally efficient for low-frequency data which are commonly installed in many smart homes. The additional penalty terms promoting sparsity and the appliance-specific constraints narrow down the feasible region, improve the algorithm performance and reduce the computational burden. We also have proposed a mixed-integer formulation for automatic state detection, and a post-processing technique to include state dynamics. 
Our algorithm effectively disambiguates appliances with similar operating modes, ensures reliable parameters estimation, and leads to accurate results, even when relatively small training data at a low sampling rate are considered.  When considering very large samples with resolution in the order of few seconds, specialized solvers are needed. One possible lead of research may be to relax the binary constraints in the BQP formulation and solve the resulting convex relaxation by using fast algorithms such as ADMM \cite{boyd2011distributed}. However, we leave investigating convex relaxations of our NILM algorithm for future work.

\section*{Acknowledgments}
We would like to thank the reviewers for their thoughtful comments that greatly helped to improve the manuscript.


\bibliographystyle{IEEEtran}
%

\bibliography{bibliography.bib}

\end{document}